\documentclass[twocolumn]{autart.RG}

\usepackage{cite}
\usepackage{amsmath}
\usepackage{amsfonts}
\usepackage{amssymb}
\usepackage{verbatim}
\usepackage{latexsym}
\usepackage{graphicx} 
\graphicspath{{./figs/}}

\newtheorem{draft-note}{Draft Note}

\newtheorem{theorem}{Theorem}
\newtheorem{proposition}{Proposition}

\newtheorem{assumption}{Assumption}
\newtheorem{remark}{Remark}

\newcommand{\mc}{\mathcal}
\newcommand{\bb}{\left(\cdot\right)}
\newcommand{\bbb}{\left(\cdot,\cdot\right)}
\newcommand{\eqbyd}{:=}
\newcommand{\C}[1]{\mathbb{#1}}
\newcommand{\N}{\C{N}}
\newcommand{\R}{\C{R}}

\newcommand{\Rnm}[2]{\R^{{#1}\times{#2}}}
\newcommand{\support}[2]{\operatorname{h}(#1,#2)}
\newcommand{\supportbb}[1]{\operatorname{h}(#1,\cdot)}
\newcommand{\gauge}[2]{\operatorname{g}(#1,#2)}
\newcommand{\gaugebb}[1]{\operatorname{g}(#1,\cdot)}

\newcommand{\HDB}[2]{\operatorname{H}_{\mc{B}^n}(#1,#2)}
\newcommand{\mcB}{\mathcal{B}}
\newcommand{\mcC}{\mathcal{C}}
\newcommand{\mcD}{\mathcal{D}}
\newcommand{\mcI}{\mathcal{I}}
\newcommand{\mcJ}{\mathcal{J}}
\newcommand{\mcK}{\mathcal{K}}
\newcommand{\mcL}{\mathcal{L}}
\newcommand{\mcP}{\mathcal{P}}
\newcommand{\mcQ}{\mathcal{Q}}
\newcommand{\mcR}{\mathcal{R}}
\newcommand{\mcS}{\mathcal{S}}
\newcommand{\mcT}{\mathcal{T}}
\newcommand{\mcX}{\mathcal{X}}
\newcommand{\mcY}{\mathcal{Y}}

\newcommand{\uP}{\underline{\mathcal{P}}}
\newcommand{\oP}{\overline{\mathcal{P}}}

\newcommand{\uC}{\underline{\mathcal{C}}}
\newcommand{\oC}{\overline{\mathcal{C}}}
\newcommand{\uQ}{\underline{\mathcal{Q}}}
\newcommand{\oQ}{\overline{\mathcal{Q}}}
\newcommand{\uV}{\underline{V}}
\newcommand{\oV}{\overline{V}}

\newcommand{\uu}{\underline{u}}
\newcommand{\ou}{\overline{u}}
\newcommand{\uz}{\underline{z}}
\newcommand{\oz}{\overline{z}}
\newcommand{\uv}{\underline{v}}
\newcommand{\ov}{\overline{v}}

\numberwithin{equation}{section}
\date{\today}
\begin{document}
\begin{frontmatter}
\runtitle{The Minkowski--Bellman Equation}
  \title{The Minkowski--Bellman Equation}
  \author[svr]{Sa\v{s}a~V.~Rakovi\'{c}\thanksref{corr}} 
	\address[svr]{Beijing Institute of Technology, Beijing, China.}
  \thanks[corr]{E--mail: sasa.v.rakovic@gmail.com. Tel.: +44 7799775366.}
\begin{abstract}
This manuscript studies the Minkowski--Bellman equation, which is the Bellman equation arising from finite or infinite horizon optimal control of unconstrained linear discrete time systems with stage and terminal cost functions specified as Minkowski functions of proper $C$--sets. In regards to the finite horizon optimal control, it is established that, under natural conditions, the Minkowski--Bellman equation and its iteration are well posed. The characterization of the value functions and optimizer maps is derived. In regards to the infinite horizon optimal control, it is demonstrated that, under the same natural conditions, the fixed point of the Minkowski--Bellman equation is unique, in terms of the value function, over the space of Minkowski functions of proper $C$--sets. The characterization of the fixed point value function and optimizer map is reported. 
\end{abstract}
\end{frontmatter}

\section{Introduction}
\label{sec:01}

Dynamic programming~\cite{bellman:1957,bellman:dreyfus:1962,bertsekas:2005,bertsekas:2018} is an indispensable mathematical technique for closed loop characterization of optimal control. The closed loop solution to finite horizon optimal control of unconstrained discrete time systems, induced by a state transition map $(x,u)\mapsto f(x,u)$, with stage and terminal cost functions,  $(x,u)\mapsto \ell(x,u)$ and $x\mapsto V_f(x)$,  can be obtained by iterating the Bellman equation~\cite{bellman:1957,bellman:dreyfus:1962} given, for all integers $k$ over the considered finite horizon and all states $x$, by
\begin{align*}
V_{k+1}(x)&=\min_{u}\ \ell(x,u)+V_k(f(x,u)),\text{ and}\\
u_{k+1}(x)&=\arg\min_{u}\ \ell(x,u)+V_k(f(x,u)),
\end{align*}
with boundary condition $V_0(x)=V_f(x)$ for all states $x$. Likewise, the properties of the fixed point of the Bellman equation~\cite{bertsekas:2005,bertsekas:2018} taking the form, for all states $x$,
\begin{align*}
V(x)&=\min_{u}\ \ell(x,u)+V(f(x,u)),\text{ and}\\
u(x)&=\arg\min_{u}\ \ell(x,u)+V(f(x,u)),
\end{align*}
play a key role in deriving the closed loop solution to the related infinite horizon optimal control problem. 

A celebrated optimal control problem that admits an elegant and easily computable solution is the linear quadratic regulator~\cite{kalman:1960,anderson:moore:1990}. Finite horizon, discrete time, linear quadratic regulator refers to finite horizon optimal control of an unconstrained linear discrete time system, $x^+=Ax+Bu$,  with stage and terminal cost functions specified as (strictly) convex in $x$ and strictly convex in $u$ quadratic functions, $\ell(x,u)=x^TQx+u^TRu$ and $V_f(x)=x^TQ_fx$. In this setting, dynamic programming produces sequences of (strictly) convex in $x$ quadratic value functions $V_{k+1}\bb$, with values $V_{k+1}(x)=x^TP_{k+1}x$, and linear optimizer functions $u_{k+1}\bb$, with values $u_{k+1}(x)=K_{k+1}x$. The corresponding value and optimizer functions are characterized by the dynamic Riccati equations specified, for all integers $k$ over the considered finite horizon, by
\begin{align*}
P_{k+1}&=Q+A^TP_kA-A^TP_kB(R+B^TP_kB)^{-1}B^TP_kA\\
K_{k+1}&=-(R+B^TP_kB)^{-1}B^TP_kA
\end{align*}
with boundary condition $P_0=Q_f$.  In the case of the infinite horizon, discrete time, linear quadratic regulator, the solution in terms of (strictly) convex in $x$ quadratic value function $V\bb$, with values $V(x)=x^TP x$, and linear optimizer function $u \bb$, with values $u(x)=K x$, is entirely determined by the algebraic Riccati equation 
\begin{align*}
P&=Q+A^TP A-A^TP B(R+B^TP B)^{-1}B^TP A\\
K&=-(R+B^TP B)^{-1}B^TP A.
\end{align*}
The detailed study, properties and practical relevance, of the finite and infinite  horizon, discrete time, linear quadratic regulator can be found in numerous references including early, but fundamental, references~\cite{kalman:1960,anderson:moore:1990}.

Peculiarly enough, even in the case of the unconstrained linear discrete time systems, the characterization of optimal control with nonquadratic stage and terminal cost functions is considerably less understood. Relevant instances of uncharted optimal control problems are the finite and infinite horizon optimal control of  unconstrained  linear discrete time systems with stage and terminal cost functions specified as Minkowski functions of proper $C$--sets. This class of optimal control problems encapsulates optimal control problems of unconstrained linear discrete time systems with stage and terminal cost functions specified as vector norms, since Minkowski functions of proper $C$--sets are nonnegative, finite valued, continuous and sublinear functions~\cite{rockafellar:1970,schneider:1993}, and vector norms can be represented via Minkowski functions of suitably defined symmetric proper $C$--sets. In these important instances, the characterization of closed loop solutions is not available in the literature. This manuscript provides characterization and computation of the closed loop solutions to both, finite and infinite horizon,  optimal control problems of unconstrained linear discrete time systems with stage and terminal cost functions specified as Minkowski functions of proper $C$--sets.  The developed results are novel and deliver a missing and relevant analogue to the celebrated linear quadratic regulator. For obvious reasons, the derived solutions can be termed as the, finite and infinite horizon, linear Minkowski regulator. The developed solution methodology enriches engineering utility of unconstrained optimal control, and also provides beneficial results within the context of \emph{inter alia} synthesis and analysis of constrained control~\cite{blanchini:miani:2008}, stabilizing control~\cite{artstein:1983,sontag:1998} and model predictive control~\cite{mayne:rawlings:rao:scokaert:2000,rawlings:mayne:2009,mayne:2014}. 

The underlying objects of study in this manuscript are the Minkowski--Bellman equation, its iteration and its fixed point. The term Minkowski--Bellman equation refers to the Bellman equation associated with the optimal control of unconstrained linear discrete time systems with stage and terminal cost functions specified as Minkowski functions of proper $C$--sets. Such a Bellman equation, its iteration and its fixed point can be studied with techniques from nonsmooth analysis~\cite{clarke:ledayaev:stern:wolenski:1998} and variational analysis~\cite{rockafellar:wets:2009}. This manuscript resorts to set--valued tools from classical results~\cite{artstein:1974,artstein:1995}. In this sense, the analysis in this manuscript makes use of the set--dynamics approach; This approach has already provided convinient tools for studies of minimality of invariant sets~\cite{artstein:rakovic:2008}, set invariance under output feedback~\cite{artstein:rakovic:2011}, and the Minkowski--Lyapunov equation~\cite{rakovic:lazar:2014,rakovic:2017}. This manuscript first analyses one step of the underlying dynamic programming iteration, namely the Minkowski--Bellman equation, and it establishes that the value function is Minkowski function of a proper $C$--set, and that its optimizer map is positively homogeneous of the first degree, compact--, convex--valued, locally bounded and outer semicontinuous map when it is set--valued, which is positively homogeneous of the first degree and continuous function when it is single--valued. The characterization of the generator set of the value function is obtained as one step iterate of a suitably defined set--dynamics. The characterization of the optimizer map is also entirely determined by this set--dynamics. 
The underlying dynamic programming iteration is in one--to--one correspondence with the iteration of the introduced set--dynamics of the generator sets. This set--dynamics is consequently utilized to characterize the iterates of the Minkowski--Bellman equation. Finally, the fixed point of the set--dynamics of the generator sets is utilized to show that the fixed point of the Minkowski--Bellman equation is unique, in terms of the value function, over the space of Minkowski functions of proper $C$--sets as well as to characterize the fixed point value function and optimizer map. 

\textbf{Manuscript Structure:}   Section~\ref{sec:02} formulates the Minkowski--Bellman equation and specifies the objectives of this manuscript. Section~\ref{sec:03} provides technical background and studies prototype problem underpinning the Minkowski--Bellman equation. Section~\ref{sec:04} studies the set--dynamics of the generator sets of the iterates of the Minkowski--Bellman equations, and examines their structural, monotonicity, boundedness, convergence and fixed point properties. Section~\ref{sec:05} considers the lower, arbitrary and upper iterates of the Minkowski--Bellman equation as well as the fixed point of the Minkowski--Bellman equation. Section~\ref{sec:06} specializes results to polytopic setting. Conclusions are drawn in Section~\ref{sec:07}.

\textbf{Typographical Convention:} We do not distinguish between a variable $x\in \R^n$ and its vectorized form. In this sense $(x,u)$ is written instead $(x^T\ u^T)^T$ and, on a few occasions, $(A\ B)(x,u)$ is written instead of $(A\ B)(x^T\ u^T)^T$.  No confusion should arise.  For clarity, some of the proofs are reported in the appendices. 

\textbf{Basic Nomenclature and Definitions:}  The set of nonnegative integers is denoted by $\N$. The set of real numbers is denoted by $\R$, while $\R_{\ge 0}$ denotes the set of nonnegative real numbers. $I$ and $0$ denote the identity and zero matrices. For $(x,u)\in\R^n\times\R^m$, the projection matrices $(x,u)\mapsto x$ and $(x,u)\mapsto u$ are
\begin{align*}
P_x&\eqbyd (I\ O)\text{ with }I\in \Rnm{n}{n}\text{ and }O\in \Rnm{n}{m},\text{ and}\\
P_u&\eqbyd (O\ I)\text{ with }O\in \Rnm{m}{n}\text{ and }I\in \Rnm{m}{m}.
\end{align*}
The Minkowski set addition of $\mcX \subseteq \R^n$ and $\mcY \subseteq \R^n$ is 
\begin{equation*}
\mcX\oplus \mcY\eqbyd \{ x+y\ :\  x\in \mc{X},\ y\in \mcY\}. 
\end{equation*} 
The image of a set $\mcX$ under a matrix (or a scalar) $M$ is
\begin{equation*}
M\mcX\eqbyd\{Mx\ :\ x\in \mcX\}.
\end{equation*}  
A set $\mcX\subseteq \R^n$ is symmetric, with respect to $0\in\R^n$, if $\mcX=-\mcX$. A set $\mcX\subset \R^n$ is \emph{a $C$--set} if it is compact, convex, and contains the origin. A set $\mcX\subset \R^n$ is \emph{a proper $C$--set} if it is a $C$--set and contains the origin in its interior. A polyhedron is the (convex) intersection of a finite number of open and/or closed half--spaces. A polytope is a closed and bounded polyhedron. For a set $\mcX$ in $\R^n$ with $0\in \mcX$,  its polar set $\mcX^*$ is given by
\begin{equation*}
\mcX^*\eqbyd \{x\in\R^n\ :\ \forall y\in\mcX,\ y^T x\le 1\}.
\end{equation*}
The support function $\supportbb{\mcX}$ of a nonempty closed convex set $\mcX$ is specified, for all $y\in \R^n$, by
\begin{equation*}
\support{\mcX}{y}\eqbyd \sup_x \{y^T x\ :\ x\in \mcX\}.
\end{equation*}
The Minkowski (gauge) function $\gaugebb{\mcX}$ of a proper $C$--set $\mcX$ in $\R^n$ is given, for all $y\in \R^n$, by
\begin{equation*}
\gauge{\mcX}{y}\eqbyd \min_\eta \{\eta\ :\ y\in \eta \mcX,\ \eta\ge 0\}.
\end{equation*}
Given any two nonempty compact subsets $\mcX$ and $\mcY$ of $\R^n$ their Hausdorff distance is 
\begin{equation*}
\HDB{\mcX}{\mcY}\eqbyd \min_{\eta\ge 0}\{ \eta\ :\  \mcX\subseteq \mcY\oplus \eta \mcB^n\text{ and }\mcY\subseteq \mcX\oplus \eta \mcB^n\}, 
\end{equation*} 
where $\mcB^n$ is the closed unit Euclidean norm ball in $\R^n$. 

A function $f\bb\ :\ \R^n\rightarrow \R^m$ is: $(a)$ positively homogeneous of the first degree if $f(\eta x)=\eta f(x)$ for all $\eta\in\R_{\ge 0}$ and all $x\in\R^n$, $(b)$ a $\pi$--class function if it is positively homogeneous of the first degree and continuous, $(c)$ subadditive if $f(x_1+x_2)\le f(x_1)+f(x_2)$ for all $x_1\in\R^n$ and $x_2\in\R^n$, and $(d)$ sublinear if it is positively homogeneous of the first degree and  subadditive.

A set--valued map $F\bb$ associates subsets $F(x)$ of $\R^m$ to points $x$ in $\R^n$. A set--valued map $F\bb\ :\ \R^n\rightrightarrows \R^m$ is:
\begin{itemize}
\item [(i)] positively homogeneous of the first degree, if $F(\eta x)=\eta F(x)$ for all $\eta\in\R_{\ge 0}$ and all $x\in \R^n$;
\item [(ii)] compact--valued at $x\in \R^n$, if $F(x)$ is a compact subset of $\R^m$;
\item [(iii)] convex--valued at $x\in \R^n$, if $F(x)$ is a convex subset of $\R^m$;
\item [(iv)] locally bounded at $x\in \R^n$, if there is a neighborhood $\mcX$ of $x$ such that the set $F(\mcX):=\cup_{x\in \mcX} F(x)$ is bounded;
\item [(v)] outer semicontinuous at $x\in \R^n$, if for every convergent sequence $x_k \rightarrow x$ and every convergent sequence $y_k\rightarrow y$ with $y_k\in F(x_k)$ it holds that $y \in F(x)$.
\end{itemize}     
A set--valued map $F\bb\ :\ \R^n\rightrightarrows \R^m$  is a $\Pi$--class set--valued map if it satisfies property $(i)$ as well as properties $(ii)$, $(iii)$, $(iv)$ and $(v)$ for all $x\in \R^n$. We note that if a $\Pi$--class set--valued map is single--valued, then it is a  $\pi$--class single--valued function by virtue of~\cite[Corollary~5.20.]{rockafellar:wets:2009}. A single--valued function $f\bb\ :\ \R^n\rightarrow \R^m$ is a selection of a set--valued map $F\bb\ :\ \R^n\rightrightarrows \R^m$ if, for all $x\in \R^n$, $f(x)\in F(x)$.

\section{Preliminaries}
\label{sec:02}

\subsection{The Linear Minkowski Regulator}
\label{sec:02.01}

The linear discrete time dynamical systems are given by
\begin{equation}
\label{eq:02.01}
x^+=Ax+Bu, 
\end{equation}
where $x\in\R^n$ and $u\in\R^m$ are the current state and current control, while $x^+\in \R^n$ is the successor state and the matrix pair $(A,B)$ is of compatible  dimensions. 
\begin{assumption}
\label{ass:02.01}
The matrix pair $(A,B)\in \Rnm{n}{n}\times \Rnm{n}{m}$ is strictly stabilizable.
\end{assumption}
Strict stabilizability signifies that the dynamics of the uncontrollable part of the system~\eqref{eq:02.01} is strictly stable. 

The stage cost function $\ell\bbb\ :\ \R^n\times\R^m\rightarrow \R_{\ge 0}$ is given, for all  $(x,u)\in\R^n\times\R^m$, by
\begin{equation}
\label{eq:02.02}
\ell(x,u)=\gauge{\mcC}{(x,u)}.
\end{equation}
\begin{assumption}
\label{ass:02.02}
The set $\mcC$ is a proper $C$--set in $\R^{n+m}$.
\end{assumption}
The terminal cost function $V_f\bb\ :\ \R^n\rightarrow \R_{\ge 0}$ is specified, for all $x\in\R^n$, by
\begin{equation}
\label{eq:02.03}
V_f(x)=\gauge{\mcQ_f}{x}.
\end{equation}

\begin{assumption}
\label{ass:02.03}
The set $\mcQ_f$ is a proper $C$--set in $\R^n$.
\end{assumption}

The finite horizon linear Minkowski regulator problem refers to the determination of a pair of finite state and control sequences, $\{x_0,x_1,\ldots,x_{N}\}$ and $\{u_0,u_1,\ldots,u_{N-1}\}$, which is dynamically consistent with the system~\eqref{eq:02.01} and an initial condition $x_0=x$ (so that $x_{k+1}=Ax_k+Bu_k,\ k\in \N_{N-1}:=\{0,1,\ldots,N-1\}$, with $x_0=x$), and which minimizes the cost function
\begin{equation}
\label{eq:02.04}
\sum_{k=0}^{N-1}\ell(x_k,u_k)+V_f(x_N).
\end{equation}

Likewise, the infinite horizon linear Minkowski regulator problem refers to the determination of a pair of infinite state and control sequences, $\{x_0,x_1,\ldots\}$ and $\{u_0,u_1,\ldots\}$, which is dynamically consistent with the system~\eqref{eq:02.01} and an initial condition $x_0=x$  (so that $x_{k+1}=Ax_k+Bu_k,\ k\in \N$, with $x_0=x$), and which minimizes the cost function
\begin{equation}
\label{eq:02.05}
\sum_{k=0}^{\infty}\ell(x_k,u_k).
\end{equation}

\subsection{The Minkowski--Bellman Equation}
\label{sec:02.02}
In analogy to the celebrated linear quadratic regulator, closed loop solutions of the finite and infinite horizon linear Minkowski regulator problems are sought. These solutions are derived by utilizing dynamic programming. The Minkowski--Bellman equation is the corresponding Bellman equation associated with optimal control of the system~\eqref{eq:02.01} with the stage and terminal cost functions represented via Minkowski functions of proper $C$--sets, as specified in~\eqref{eq:02.02} and~\eqref{eq:02.03}. The Minkowski--Bellman equation takes the form, for all $k\in\N$  and all $x\in \R^n$,
\begin{subequations}
\label{eq:02.06}
\begin{align}
\label{eq:02.06.a}
V_{k+1}(x)&\eqbyd \min_u\ \ell(x,u)+V_k(Ax+Bu),\text{ and }\\
\label{eq:02.06.b}
u_{k+1}(x)&\eqbyd \arg\min_u\ \ell(x,u)+V_k(Ax+Bu),
\end{align}
\end{subequations}
with the boundary condition given, for all $x\in\R^n$, by
\begin{equation}
\label{eq:02.07}
V_0(x)\eqbyd V_f(x).
\end{equation} 
Throughout this manuscript, $V_k\bb$ and $u_k\bb$ are referred to as the value function and its optimizer map.

 We also consider the fixed point of the above Minkowski--Bellman equation taking the form, for  all $x\in \R^n$,
\begin{subequations}
\label{eq:02.08}
\begin{align}
\label{eq:02.08.a}
V(x)&=\min_u\ \ell(x,u)+V(Ax+Bu),\text{ and}\\
\label{eq:02.08.b}
u(x)&= \arg\min_u \ell(x,u)+V(Ax+Bu),
\end{align}
\end{subequations}
where the value function $V\bb$ and its optimizer map $u\bb$ are to be determined.

\subsection{Problem Description}
\label{sec:02.03}

Our main chore is to characterize the solution, and discuss topological properties, of the Minkowski--Bellman equation, its iteration and its fixed point. 

The first goal is to establish that the value functions $V_{k}\bb$ are Minkowski functions of proper $C$--sets in $\R^n$, and that each of the optimizer maps $u_{k}\bb$ is a $\Pi$--class set--valued map. This goal also requires one to characterize the generator sets $\mcP_k$ of the value functions $V_{k}\bb$ as well as to characterize the optimizer maps $u_{k}\bb$. 

The second goal is to establish that the fixed point of the Minkowski--Bellman equation  admits a unique solution in terms of the value function $V\bb$ over the space of Minkowski functions of a proper $C$--sets in $\R^n$ as well as to show that related optimizer map $u\bb$ is a $\Pi$--class set--valued map. This goal also requires one to characterize the generator set  $\mcP$ of the value function  $V\bb$ as well as to characterize the optimizer map $u\bb$.

\section{Background and Prototype Problem}
\label{sec:03}

\subsection{Background}
\label{sec:02.05}

Our analysis uses the following fundamental results~\cite{rockafellar:1970,schneider:1993}.

\begin{theorem}~\cite[Theorem~1.6.1.]{schneider:1993}
\label{thm:03.01}
Let $\mcX$ be a proper $C$--set in $\R^n$. Then its polar set $\mcX^*$ is itself a proper $C$--set in $\R^n$ and it holds that $(\mcX^*)^*=\mcX$.
\end{theorem}
\begin{theorem}~\cite[Theorem~1.7.1.]{schneider:1993}
\label{thm:03.02}
If $f\bb\ :\ \R^n\rightarrow \R$ is a sublinear function, then there is a unique nonempty convex compact subset of $\R^n$ with support function $f\bb$.
\end{theorem}
\begin{theorem}~\cite[Theorem~1.7.6.]{schneider:1993}
\label{thm:03.03}
For a proper $C$--set $\mcX$ in $\R^n$ it holds that
\begin{equation}
\label{eq:03.01}
\forall x\in\R^n,\quad \gauge{\mcX}{x}=\support{\mcX^*}{x},
\end{equation}
where $\mcX^*$ is the polar set of $\mcX$.
\end{theorem}
The following two well--known facts also prove helpful.
\begin{proposition}
\label{prop:03.01}
Let $\mcX$ and $\mcY$ be two nonempty closed convex sets in $\R^n$, and let also $M\in\Rnm{n}{p}$. Then
\begin{subequations}
\label{eq:03.02}
\begin{align}
\label{eq:03.02.a}
\forall y\in\R^p,&\quad \support{\mcX}{My}=\support{M^T\mc{X}}{y},\text{ and}\\
\label{eq:03.02.b}
\forall z\in \R^n,&\quad \support{\mcX\oplus \mcY}{z}=\support{\mcX}{z}+\support{\mcY}{z}.
\end{align}
\end{subequations}
\end{proposition}
\begin{proposition}
\label{prop:03.02}
Let $\mc{X}$ and $\mc{Y}$ be proper $C$--sets in $\R^n$.  Then $\mc{X}\subseteq \mc{Y}$ if and only if 
\begin{equation}
\label{eq:03.03}
\forall x\in\R^n,\quad \support{\mcX}{x}\le \support{\mcY}{x},
\end{equation}
or, equivalently, if and only if
\begin{equation}
\label{eq:03.04}
\forall x\in\R^n,\quad \gauge{\mcY}{x}\le \gauge{\mcX}{x}.
\end{equation}
\end{proposition}

\subsection{Prototype Problem}
\label{sec:03.01}

The prototype problem provides the characterization of the successor value function $V_{k+1}\bb$ and its optimizer map $u_{k+1}\bb$ when the stage cost function $\ell\bbb$ and the current value function $V_k\bb$ are Minkowski functions of proper $C$--sets. The prototype problem is a parametric optimization problem, with respect to $x\in \R^n$, 
\begin{subequations}
\label{eq:03.05}
\begin{align}
\label{eq:03.05.a}
V^+(x)&\eqbyd \min_u\ \ell(x,u)+V(Ax+Bu),\text{ and }\\
\label{eq:03.05.b}
u^+(x)&\eqbyd \arg\min_u\ \ell(x,u)+V(Ax+Bu),
\end{align}
\end{subequations}
where, for all $(x,u)\in \R^n\times\R^m$ and all $x\in \R^n$, 
\begin{equation}
\label{eq:03.06}
\ell(x,u)\eqbyd\gauge{\mcC}{(x,u)}\text{ and }V(x)\eqbyd\gauge{\mcP}{x},
\end{equation}
and $\mcC$ and $\mcP$ are proper $C$--sets in $\R^{n+m}$ and $\R^n$.

The solution to the above prototype problem can be derived by employing the set--dynamics of the generator sets $\mcP^+$ of the successor value function $V^+\bb$ in terms of generator sets, $\mcC$ and $\mcP$, of the stage cost and current value functions, $\ell\bbb$ and $V\bb$, as specified by
\begin{equation}
\label{eq:03.07}
\mcP^+:=P_x\left(\mcC^*\oplus (A\ B)^T\mcP^*\right)^*.
\end{equation}
\begin{theorem}
\label{thm:03.04}
Take any $(A,B)\in\Rnm{n}{n}\times\Rnm{n}{m}$, and any proper $C$--sets $\mcC$ and $\mcP$ in $\R^{n+m}$ and $\R^n$. Let 
\begin{equation}
\mcT^+:=\left(\mcC^*\oplus (A\ B)^T\mcP^*\right)^*\text{ and } \mcP^+=P_x\mcT^+.
\label{eq:03.08}
\end{equation}
Consider the prototype problem~\eqref{eq:03.05}--\eqref{eq:03.06}.
\begin{itemize}
\item [(i)] $\mcT^+$ and $\mcP^+$ are proper $C$--sets in $\R^{n+m}$ and $\R^{n}$. 
\item [(ii)] $V^+\bb\ :\ \R^n\rightarrow \R_{\ge 0}$ is the Minkowski function of proper $C$--set $\mcP^+$ so that, for all $x\in\R^n$,
\begin{equation}
\label{eq:03.09}
V^+(x)=\gauge{\mcP^+}{x}.
\end{equation}	
\item [(iii)] $u^+\bb\ :\ \R^n\rightrightarrows \R^m$ is a $\Pi$--class set--valued map equivalently given, for all $x\in \R^n$, by\footnote{For proper $C$--sets $\mcP^+$ and $\mcT^+$ in $\R^n$ and $\R^{n+m}$, and for $(x,u)\in\R^n\times \R^m$, we frequently write $(x,u)\in\gauge{\mcP^+}{x}\mcT^+$, which is equivalent to $\frac{1}{\gauge{\mcP^+}{x}}(x,u)\in \mcT^+$ for $x\neq 0$ and $(0,u)\in \{(0,0)\}$ (i.e. $u=0$) for $x=0$.}
\begin{equation}
\label{eq:03.10}
u^+(x)=\{u\in\R^m\ :\ (x,u)\in \gauge{\mcP^+}{x}\mcT^+\}.
\end{equation}	
\item [(iv)] $\nu^+\bb\ :\ \R^n\rightarrow \R^m$ given, for all $x\in \R^n$, by
\begin{equation}
\label{eq:03.11}
\nu^+(x)=\arg\min_{u}\{u^Tu\ :\ u\in u^+(x)\}
\end{equation}
is a $\pi$--class function such that, for all $x\in\R^n$,
\begin{equation}
\label{eq:03.12}
\nu^+(x)\in u^+(x).
\end{equation}
\end{itemize}
\end{theorem}
Theorem~\ref{thm:03.04} implies directly that the Minkowski--Bellman equation~\eqref{eq:02.06} with the boundary condition~\eqref{eq:02.07} generates sequence of the value functions $V_{k+1}\bb$, terms of which are Minkowski functions of proper $C$--sets in $\R^n$; It also generates a sequence of related optimizer maps $u_{k+1}\bb$, terms of which are $\Pi$--class set--valued maps that admit $\pi$--class selections $\nu_{k+1}\bb$.
\begin{remark}
\label{rem:03.01}
The Minkowski function of a proper $C$--set is $0$ for $x=0$ and striclty positive for all $x\in \R^n\setminus \{0\}$, finite valued, continuous and sublinear function~\cite{rockafellar:1970,schneider:1993}. A $\Pi$--class set--valued map, which is single--valued is a $\pi$--class (i.e. positively homogeneous of the first degree and continuous) function. These generic properties of the successor value function $V^+\bb$, its optimizer map $u^+\bb$ and selection $\nu^+\bb$ as well as of the value functions $V_{k+1}\bb$, their optimizer maps $u_{k+1}\bb$ and  selections $\nu_{k+1}\bb$ are granted by definitions. In what follows, unless necessary, these inherent topological properties are neither formally stated nor elaborated on. 
\end{remark}

\section{Set--Dynamics of Generator Sets}
\label{sec:04}

With Theorem~\ref{thm:03.04} in mind, the iterates of the Minkowski--Bellman equation~\eqref{eq:02.06} with the boundary condition~\eqref{eq:02.07} are determined by the iterates of the set--dynamics~\eqref{eq:03.07} taking the form, for all $k\in \N$,
\begin{equation}
\label{eq:04.01}
\mcP_{k+1}=P_x\left(\mcC^*\oplus (A\ B)^T\mcP_k^*\right)^*,
\end{equation}
with the initial condition $\mcP_0\eqbyd \mcQ_f$ in its polar form
\begin{equation}
\label{eq:04.02}
\mcP_0^*\eqbyd \mcQ_f^*.
\end{equation}
By the same token, the properties of the fixed point of the Minkowski--Bellman equation~\eqref{eq:02.08} are determined by the properties of the fixed point of the set--dynamics~\eqref{eq:03.07}, namely the fixed point set--equation
\begin{equation}
\label{eq:04.03}
\mcP=P_x\left(\mcC^*\oplus (A\ B)^T\mcP^*\right)^*.
\end{equation}
In view of this one--to--one correspondence between the set--dynamics~\eqref{eq:03.07}, its iteration~\eqref{eq:04.01}--\eqref{eq:04.02}, and its fixed point~\eqref{eq:04.03} on one hand and the Minkowski--Bellman equation (captured by the prototype problem~\eqref{eq:03.05}--\eqref{eq:03.06}), its iteration~\eqref{eq:02.06}--\eqref{eq:02.07}, and its fixed point~\eqref{eq:02.08} on the other hand, we discuss key properties of the iterates of the set--dynamics~\eqref{eq:04.01}, and subsequently translate these properties to the corresponding iterates of the Minkowski--Bellman equation. The iterates of the set--dynamics~\eqref{eq:04.01} are examined for arbitrary initial conditions $\mcP_0$ specified via related polar sets $\mcP_0^*$. 

The following preliminary result proves very useful. 
\begin{proposition}
\label{prop:04.01}
Suppose Assumptions~\ref{ass:02.01} and ~\ref{ass:02.02} hold.
\begin{itemize}
\item [(i)] There exists a proper $C$--set $\mcL$ in $\R^n$, for which $\gauge{\mcL}{\cdot}$ verifies control Lyapunov decrease condition, i.e. for all $x\in\R^n$, there exists a $u\in\R^m$ such that
\begin{equation}
\label{eq:04.04}
\gauge{\mcL}{Ax+Bu}+\gauge{\mcC}{(x,u)}\le \gauge{\mcL}{x}.
\end{equation}
\item [(ii)] If a proper $C$--set $\mcL$ in $\R^n$ verifies~\eqref{eq:04.04}, and if  
\begin{align}
\mcK^+=\left(\mcC^*\oplus (A\ B)^T\mcL^*\right)^*\text{ and } \mcL^+&=P_x\mcK^+,
\label{eq:04.05}
\end{align}
then $\mcK^+$ and $\mcL^+$ are proper $C$--sets in $\R^{n+m}$ and $\R^n$, respectively and,  for all $(x,u)\in\R^n\times\R^m$,
\begin{equation}
\label{eq:04.06}
\gauge{\mcK^+}{(x,u)}=\gauge{\mcL}{Ax+Bu}+\gauge{\mcC}{(x,u)},
\end{equation}
and, for all $x\in\R^n$, there exists a $u\in\R^m$ such that
\begin{equation}
\label{eq:04.07}
\gauge{\mcL^+}{x}\le \gauge{\mcK^+}{(x,u)}\le \gauge{\mcL}{x}.
\end{equation}
\item [(ii)] If Assumption~\ref{ass:02.03} also holds, then the proper $C$--set $\mcL$ in $\R^n$ from the above assertion $(i)$ can be chosen to satisfy
\begin{equation}
\label{eq:04.08}
\mcL\subseteq \mcQ_f.
\end{equation}
\end{itemize}
\end{proposition}

\subsection{Characterization, Monotonicity and Boundedness}
\label{sec:04.01}

Theorem~\ref{thm:03.04} and mathematical induction reveal the structure of the generator sets $\mcP_k$ as well as related sets $\mcT_k$. 
\begin{proposition}
\label{prop:04.02}
Suppose Assumptions~\ref{ass:02.01} and~\ref{ass:02.02} hold. Let, for all $k\in \N$,
\begin{equation}
\label{eq:04.09}
\mcT_{k+1}:=\left(\mcC^*\oplus (A\ B)^T\mcP_k^*\right)^*.
\end{equation}
If $\mcP_0^*$ is a $C$--set in $\R^n$, then, for all $k\in \N$, sets $\mcT_{k+1}$ and $\mcP_{k+1}$ are proper $C$--set in $\R^{n+m}$ and $\R^n$.
\end{proposition}
As long as $\mcP_0^*$ is a $C$--set in $\R^n$, Proposition~\ref{prop:04.02} and Theorem~\ref{thm:03.04} characterize completely the value functions $V_k\bb$, optimizer maps $u_k\bb$ and related selections $\nu_k\bb$. The value functions $V_k\bb$ are Minkowski functions of proper $C$--sets $\mcP_k$, the optimizer maps $u_k\bb$ are $\Pi$--class set--valued maps and the related selections $\nu_k\bb$ are $\pi$--class functions. In particular, for all $k\in\N$ and all $x\in \R^n$,
\begin{align}
V_{k+1}(x)&=\gauge{\mcP_{k+1}}{x},\nonumber\\
u_{k+1}(x)&=\{u\in\R^m\ :\ (x,u)\in \gauge{\mcP_{k+1}}{x}\mcT_{k+1}\},\text{ and}\nonumber\\
\label{eq:04.10}
\nu_{k+1}(x)&=\arg\min_u\{u^Tu\ :\ u\in u_{k+1}(x)\}.
\end{align}

The iterates of the set--dynamics~\eqref{eq:04.01} preserve monotonicity, as established by the following result.
\begin{proposition}
\label{prop:04.03}
Suppose Assumptions~\ref{ass:02.01} and~\ref{ass:02.02} hold. 
\begin{itemize}
\item [(i)] If  $\mcP_0^*$ is a $C$--set in $\R^n$ such that $\mcP_0^*\subseteq \mcP_1^*$, then, for all $k\in \N$, 
\begin{equation}
\label{eq:04.11}
\mcP_{k+2}\subseteq \mcP_{k+1}.
\end{equation}
\item [(ii)] If  $\mcP_0^*$ is a $C$--set in $\R^n$ such that $\mcP_1^*\subseteq \mcP_0^*$, then, for all $k\in \N$,
\begin{equation}
\label{eq:04.12}
\mcP_{k+1}\subseteq \mcP_{k+2}.
\end{equation}
\end{itemize}
\end{proposition}
By Propositions~\ref{prop:03.02} and~\ref{prop:04.03}, when $\mcP_0^*$ is a $C$--set in $\R^n$ such that $\mcP_0^*\subseteq \mcP_1^*$, the sequence of the value functions  $V_k\bb$ is monotonically nondecreasing with respect to $k$, i.e., for all $k\in\N$ and all $x\in\R^n$, 
\begin{equation}
\label{eq:04.13}
V_{k+1}(x)\le V_{k+2}(x).
\end{equation}
Likewise, when $\mcP_0^*$ is a $C$--set in $\R^n$ such that $\mcP_1^*\subseteq \mcP_0^*$, the sequence of the value functions $V_k\bb$ is monotonically nonincreasing with respect to $k$, i.e., for all $k\in\N$ and all $x\in\R^n$, 
\begin{equation}
\label{eq:04.14}
V_{k+2}(x)\le V_{k+1}(x).
\end{equation}
\begin{remark}
\label{rem:04.01}
When $\mcP_0^*$ is a $C$--set in $\R^n$ such that $\mcP_0^*\subseteq \mcP_1^*$, the value functions $V_k\bb$ satisfy, for all $k\in \N$, all $x\in \R^n$ and all $u\in u_{k+1}(x)$ (including $u=\nu_{k+1}(x)$), 
\begin{equation}
\label{eq:04.15}
V_{k}(x) \le \ell(x,u)+V_k(Ax+Bu)=V_{k+1}(x).
\end{equation}
Similarly, when $\mcP_0^*$ is a $C$--set in $\R^n$ such that $\mcP_1^*\subseteq \mcP_0^*$, the value functions $V_k\bb$ preserve the \emph{strong Lyapunov decrease} property. Namely, for all $k\in \N$, all $x\in \R^n$ and all $u\in u_{k+1}(x)$ (including $u=\nu_{k+1}(x)$), 
\begin{equation}
\label{eq:04.16}
V_{k+1}(x)=V_k(Ax+Bu)+\ell(x,u)\le V_{k}(x).
\end{equation}
\end{remark}

The iterates of the set--dynamics~\eqref{eq:04.01} are also monotone with respect to their initial conditions.
\begin{proposition}
\label{prop:04.04}
Suppose Assumptions~\ref{ass:02.01} and~\ref{ass:02.02} hold. If  $(\mcP_0')^*$ and $(\mcP_0'')^*$ are $C$--sets in $\R^n$ such that $(\mcP_0'')^*\subseteq (\mcP_0')^*$, then, for all $k\in \N$, 
\begin{equation}
\label{eq:04.17}
\mcP'_{k+1}\subseteq \mcP''_{k+1}.
\end{equation}
\end{proposition}

The iterates of the set--dynamics~\eqref{eq:04.01} are also suitably inner and outer bounded by constant proper $C$--sets.
\begin{proposition}
\label{prop:04.05}
Suppose Assumptions~\ref{ass:02.01} and~\ref{ass:02.02} hold. Suppose also that $\mcL$ is a proper $C$--set in $\R^n$ that verifies relation~\eqref{eq:04.04}. Let $\mcD^+\eqbyd P_x\mcC$. If $\mcP_0^*$ is a $C$--set in $\R^n$ such that $\mcP_0^*\subseteq \mcL^*$, then, for all $k\in \N$, 
\begin{equation}
\label{eq:04.18}
\mcL\subseteq \mcP_{k+1}\subseteq \mcD^+.
\end{equation}
\end{proposition}
By Propositions~\ref{prop:03.02} and~\ref{prop:04.05}, for all $k\in\N$ and all $x\in\R^n$, 
\begin{equation}
\label{eq:04.19}
\gauge{\mcD^+}{x}\le V_{k+1}(x)\le \gauge{\mcL}{x}.
\end{equation}

\subsection{Convergence and Fixed Point}
\label{sec:04.02}

The convergence of the sequences of the value functions $V_k\bb$, optimizer maps $u_k\bb$ and  selections $\nu_k\bb$ is dictated by the convergence of the sequences of the generator sets $\mcP_k$ and proper $C$--sets $\mcT_k$ defined in~\eqref{eq:04.09}. 
\begin{theorem}
\label{thm:04.01}
Suppose Assumptions~\ref{ass:02.01} and~\ref{ass:02.02} hold.  If $\mcP_0^*$ is a $C$--set in $\R^n$ such that either $\mcP_0^*\subseteq \mcP_1^*$ or $\mcP_1^*\subseteq \mcP_0^*$, then
\begin{itemize}
\item [(i)] The sequence $\{\mcP_k\}_{k\ge 1}$ of the generator sets $\mcP_k$ converges\footnote{In this manuscript, the convergence of sequences of nonempty convex compact subsets of $\R^n$ is considered with respect to the Hausdorff distance.} to a proper $C$--set $\mcP_\infty$ in $\R^n$.
\item [(ii)] The sequence $\{\mcT_k\}_{k\ge 1}$ of the sets $\mcT_k$ generated by~\eqref{eq:04.09} converges to a proper $C$--set $\mcT_\infty$ in $\R^{n+m}$ satisfying
\begin{equation}
\label{eq:04.20}
\mcT_\infty=\left(\mcC^*\oplus (A\ B)^T\mcP_\infty^*\right)^*\text{ and }\mcP_\infty=P_x\mcT_\infty.
\end{equation}
\item [(iii)] The limit $\mcP_\infty$ solves the fixed point set--equation~\eqref{eq:04.03}.
\end{itemize} 
\end{theorem}

The convergence of the sequence of the generator sets $\mcP_k$ yields directly uniform convergence of the value functions $V_k\bb$ on the unit sphere $\mathbb{S}^{n-1}$. As is customary, $\mathbb{S}^{n-1}:=\{x\in\R^n\ :\ x^Tx=1\}$ denotes the unit sphere of the Euclidean  norm in $\R^n$. When $\mcP_0^*$ is a $C$--set in $\R^n$ such that either $\mcP_0^*\subseteq \mcP_1^*$ or $\mcP_1^*\subseteq \mcP_0^*$, the sequence $\{V_{k+1}\bb\}_{k\in\N}$ of the value functions converges uniformly on the unit sphere $\mathbb{S}^{n-1}$ to the Minkowski function of a proper $C$--set $\mcP_\infty$ in $\R^n$
\begin{equation*}
V_k\bb \rightarrow \gaugebb{\mcP_\infty}\text{ uniformly on }\mathbb{S}^{n-1}\text{ as }k\rightarrow \infty,
\end{equation*}
where $\mcP_\infty$ is the limit established in Theorem~\ref{thm:04.01}$(i)$. The convergence of the sequences of the generator sets $\mcP_k$ and related sets $\mcT_k$ results in uniform convergence of the optimizer maps $u_k\bb$ and the related selections $\nu_k\bb$ over the unit sphere $\mathbb{S}^{n-1}$. The limiting optimizer map $u_\infty\bb$ is a $\Pi$--class set--valued map and the related limiting selection $\nu_\infty\bb$ is a $\pi$--class function characterized by the limiting generator set $\mcP_\infty$ and the related limiting set $\mcT_\infty$. In particular, $u_\infty\bb$ and $\nu_\infty\bb$ are given, for all $x\in\R^n$, by
\begin{align}
u_\infty(x)&=\{u\in\R^m\ :\ (x,u)\in \gauge{\mcP_\infty}{x}\mcT_\infty\},\text{ and}\nonumber\\
\label{eq:04.21}
\nu_\infty(x)&=\arg\min_u\{u^Tu\ :\ u\in u_\infty(x)\}.
\end{align}
When $\mcP_0^*$ is a $C$--set in $\R^n$ such that either $\mcP_0^*\subseteq \mcP_1^*$ or $\mcP_1^*\subseteq \mcP_0^*$ the sequences $\{u_{k+1}\bbb\}_{k\in\N}$ and $\{\nu_{k+1}\bbb\}_{k\in\N}$ of the optimizer maps and related selections converge uniformly on unit sphere $\mathbb{S}^{n-1}$ to the limiting optimizer map and selection $u_\infty\bb$ and $\nu_\infty\bb$:
\begin{align*}
u_k\bb &\rightarrow u_\infty\bb\text{ uniformly on }\mathbb{S}^{n-1}\text{ as }k\rightarrow \infty\text{ and}\\
\nu_k\bb &\rightarrow \nu_\infty\bb\text{ uniformly on }\mathbb{S}^{n-1}\text{ as }k\rightarrow \infty.
\end{align*}
The technical details justifying the above convergence conclusion are provided in Appendices~C and~D.

It worth noting that, due to positive homogeneity of the first degree of the value functions,  optimizer maps and related selections, $V_k\bb$, $u_{k}\bb$ and $\nu_{k}\bb$, their uniform convergence over the unit sphere $\mathbb{S}^{n-1}$ implies directly their pointwise convergence over $\R^n$ as well as uniform convergence over nonempty compact subsets of $\R^n$. The limiting value function, optimizer map and related selction are $\gaugebb{\mcP_\infty}$, $u_\infty\bb$ and $\nu_\infty\bb$.

\begin{remark}
\label{rem:04.02}
When $\mcP_0^*$ is a $C$--set in $\R^n$ such that either $\mcP_0^*\subseteq \mcP_1^*$ or $\mcP_1^*\subseteq \mcP_0^*$, due to Theorems~\ref{thm:03.04} and~\ref{thm:04.01}$(iii)$, the limiting value function and optimizer map (or related selection), $\gaugebb{\mcP_\infty}$ and $u_\infty\bb$ (or $\nu_\infty\bb$),  form a solution to the fixed point of the Minkowski--Bellman equation~\eqref{eq:02.08}. The limiting value function $\gaugebb{\mcP_\infty}$ and its optimizer map $u_\infty\bb$ also ensure a \emph{strong Lyapunov decrease property}, i.e., for all $x\in \R^n$ and all $u\in u_\infty(x)$ (including $u=\nu_\infty(x)$), 
\begin{equation}
\label{eq:04.22}
\gauge{\mcP_\infty}{x}=\gauge{\mcC}{(x,u)}+\gauge{\mcP_\infty}{Ax+Bu}.
\end{equation} 
\end{remark}

The preceding analysis has established that, under the related monotonicity hypothesis, the sequence of the generator sets $\mcP_k$ is convergent, and that its limit is a proper $C$--set $\mcP_\infty$ in $\R^n$, which, in addition, solves the fixed point set--equation~\eqref{eq:04.03}. However, it has not been established that the limit $\mcP_\infty$ is independent of the initial conditions $\mcP_0$. That particular aspects is addressed in Section~\ref{sec:05}. The necessary analysis benefits from the fact that  taking limits does not destroy monotonicity as formally stated by the following result, which follows directly from Proposition~\ref{prop:04.04} and Theorem~\ref{thm:04.01}.
\begin{theorem}
\label{thm:04.02}
Suppose Assumptions~\ref{ass:02.01} and~\ref{ass:02.02} hold.  Suppose also that $(\mcP_0')^*$ and $(\mcP_0'')^*$ are $C$--set in $\R^n$ such that
\begin{itemize}
\item [(I)] $(\mcP_0'')^*\subseteq (\mcP_0')^*$.
\item [(II)]  $(\mcP_0')^*\subseteq (\mcP_1')^*$ or $(\mcP_1')^*\subseteq (\mcP_0')^*$.
\item [(III)] $(\mcP_0'')^*\subseteq (\mcP_1'')^*$ or $(\mcP_1'')^*\subseteq (\mcP_0'')^*$.
\end{itemize}
For the limits $\mcP_\infty'$ and $\mcP_\infty''$ established in Theorem~\ref{thm:04.01}
\begin{equation}
\label{eq:04.23}
\mcP_\infty'\subseteq \mcP_\infty''.
\end{equation}
\end{theorem}

\subsection{Symmetry of Iterates and Finite Determination}
\label{sec:04.03}

It is of interest to comment on symmetry of the iterates of the Minkowski--Bellman equation and its fixed point. 
\begin{remark}
\label{rem:04.03}
A sufficient, but not necessary, condition for symmetry of the generator sets $\mcP_k$ and limiting generator set $\mcP_\infty$ is the requirement for symmetry of the generator set $\mcC$ of the stage cost function $\ell\bbb$ and the initial condition $\mcP_0$ (or related polar set $\mcP_0^*$). 
\end{remark}

The finite determination of the fixed point set $\mcP_\infty$ is of a lot of interest for  structural reasons. 
\begin{remark}
\label{rem:04.04}
The limit $\mcP_\infty$ is finitely determined if and only if for some finite integer $k^*$ it holds that
\begin{equation}
\label{eq:04.24}
\mcP_{k^*+1}=\mcP_{k^*+2},
\end{equation}
in which case
\begin{equation}
\label{eq:04.25}
\mcP_\infty=\mcP_{k^*+2}.
\end{equation}
\end{remark}

\section{Iterates of the Minkowski--Bellman Equation}
\label{sec:05}

\subsection{Consistently Improving Lower Iterates}
\label{sec:05.01}

The iterates of the Minkowski--Bellman equation~\eqref{eq:02.06} with the boundary condition given, for all $x\in\R^n$, by 
\begin{equation}
\label{eq:05.01}
V_0(x)\eqbyd 0,
\end{equation}
are referred to as the lower iterates. The solution of the Minkowski--Bellman equation~\eqref{eq:02.06} with the boundary condition~\eqref{eq:05.01} is entirely determined by the associated generator sets produced by the set--dynamics~\eqref{eq:04.01} with the initial condition in its polar form
\begin{equation}
\label{eq:05.02}
\mcP_0^*\eqbyd\{0\}.
\end{equation}

The generator sets $\mcP_k$ have properties established in Propositions~\ref{prop:04.02},~\ref{prop:04.03}$(i)$,~\ref{prop:04.05} and Theorem~\ref{thm:04.01}. In this sense, the sequence of the generators sets $\mcP_k$ is sequence of proper $C$--sets in $\R^n$. The generator sets are inner and outer bounded by proper $C$--sets $\mcL$ and $\mcD^+$ in $\R^n$, as specified in~\eqref{eq:04.18}. The sequence of the generator sets $\mcP_k$ is also sequence of monotonically nonincreasing sets, as specified in~\eqref{eq:04.11}. This sequence is, therefore, convergent and the limit $\mcP_\infty$ of the generator sets $\mcP_k$ is proper $C$--sets in $\R^n$. The limiting generator set $\mcP_\infty$ is a solution to the fixed point set--equation~\eqref{eq:04.03}.

The value functions $V_k\bb$ are Minkowski functions of proper $C$--sets $\mcP_k$ in  $\R^n$, as specified by~\eqref{eq:04.10}. The optimizer maps $u_k\bb$ are $\Pi$--class set--valued maps and the related selections $\nu_k\bb$ are $\pi$--class functions; The optimizer maps $u_k\bb$ and the related selections $\nu_k\bb$ are characterized in~\eqref{eq:04.10}. The sequence of the value functions  are monotonically nondecreasing with respect to $k$, as specified in~\eqref{eq:04.13}. The value functions $V_k\bb$ are also lower and upper bounded by Minkowski functions of proper $C$--sets $\mcD^+$ and $\mcL$, as specified in~\eqref{eq:04.19}. The sequence of the value functions $V_k\bb$ converges uniformly over the unit sphere $\mathbb{S}^{n-1}$ to the Minkowski function of proper $C$--set $\mcP_\infty$ in $\R^n$.  The sequences of the optimizer maps $u_k\bb$ and related selections $\nu_k\bb$ converge uniformly over the unit sphere $\mathbb{S}^{n-1}$ to the limiting optimizer map and related selection $u_\infty\bb$ and $u_\infty\bb$, which are, respectively, a $\Pi$--class set--valued map and a $\pi$--class single--valued function specified by~\eqref{eq:04.21}. The limiting value function and optimizer map (or related selection) $\gaugebb{\mcP_\infty}$ and $u_\infty\bb$ (or $\nu_\infty\bb$) satisfy the fixed point of the Minkowski--Bellman equation~\eqref{eq:02.08}.

\subsection{Consistently Improving Upper Iterates}
\label{sec:05.02}

The iterates of the Minkowski--Bellman equation~\eqref{eq:02.06} with the boundary condition given, for all $x\in\R^n$, by
\begin{equation}
\label{eq:05.03}
V_0(x)\eqbyd  \gauge{\mcL}{x},
\end{equation}
are referred to as the upper iterates. The solution of the Minkowski--Bellman equation~\eqref{eq:02.06} with the boundary condition~\eqref{eq:05.03} is entirely determined by the associated generator sets produced by the set--dynamics~\eqref{eq:04.01} with the initial condition in its polar form
\begin{equation}
\label{eq:05.04}
\mcP_0^*=\mcL^*.
\end{equation}
Above, the proper $C$--set $\mcL$ in $\R^n$ satisfies relation~\eqref{eq:04.04}. 

In this case, the generator sets $\mcP_k$ have properties established in Propositions~\ref{prop:04.02},~\ref{prop:04.03}$(ii)$,~\ref{prop:04.05} and Theorem~\ref{thm:04.01}. Thus, the generator sets $\mcP_k$ possess all the properties established for the case of the lower iterates of the Minkowski--Bellman equation with one difference. Namely, in this case, the generator sets $\mcP_k$ are monotonically nondecreasing with respect to $k$, as specified in~\eqref{eq:04.12}. 

The value functions, optimizer maps and selections, $V_k\bb$, $u_k\bb$ and $\nu_k\bb$, posses all the properties established for the case of the lower iterates of the Minkowski--Bellman equation with one difference. Namely, in this setting, the value functions $V_k\bb$ are monotonically nonincreasing with respect to $k$, as specified in~\eqref{eq:04.14}. 

\subsection{Equality of Limits of Lower and Upper Iterates}
\label{sec:05.03}

As already pointed out, there is no \emph{a priori} guarantee that the limits of generator sets of the lower and upper iterates of the Minkowski--Bellman equation are equal. Hence, we proceed to demonstrate that these limits are actually identical. Let $\uP_\infty$ denote the limit of the generator sets $\uP_k$ of the lower iterates, and let $\oP_\infty$ denote the limit of the generator sets $\oP_k$ of the upper iterates.

Consider first, an  arbitrary infinite horizon optimal control process associated with the limit of the lower iterates of the Minkowski--Bellman equation. Hence, for any $x\in \R^n$, consider the infinite state and control sequences, $\mathbf{\uz}_\infty(x)=\{\uz_0(x),\uz_1(x),\ldots\}$ and $\mathbf{\uv}_\infty(x)=\{\uv_0(x),\uv_1(x),\ldots\}$, generated, for all $k\in \N$, by
\begin{subequations}
\label{eq:05.05}
\begin{align}
\label{eq:05.05.a}
\uz_{k+1}(x)&\eqbyd A\uz_{k}(x)+B\uv_{k}(x)\text{ with }\uz_{0}(x)\eqbyd x,\\
\label{eq:05.05.b}
\uv_k(x)&\in \uu_\infty(\uz_k(x)),
\end{align}
\end{subequations}
where the selection of controls $\uv_k(x)\in \uu_\infty(\uz_k(x))$ is arbitrary. In view of Remark~\ref{rem:04.02}, the sequences $\mathbf{\uz}_\infty(x)$ and $\mathbf{\uv}_\infty(x)$ converge to $(0,0)\in \R^n\times\R^m$ for any $x\in \R^n$. Hence, the limit of the sum of the values of the associated stage cost function $\ell(\uz_j(x),\uv_j(x))$ specified, for all $x\in\R^n$, by 
\begin{align}
\uV_\infty(x)&=\uV_\infty(\mathbf{\uz}_\infty(x),\mathbf{\uv}_\infty(x))\nonumber\\
\label{eq:05.06}
&\eqbyd \lim_{k\rightarrow \infty}\left(\sum_{j=0}^{k-1} \ell(\uz_j(x),\uv_j(x))\right)
\end{align}
converges pointwise in $\R^n$ as $k\rightarrow \infty$ to
\begin{equation}
\label{eq:05.07}
\uV(x)= \gauge{\uP_\infty}{x}.
\end{equation}
In fact, due to~\cite[Theorem~1.18.12.]{schneider:1993} $\uV_\infty\bb$ converges uniformly over the unit sphere $\mathbb{S}^{n-1}$ to $\gaugebb{\uP_\infty}$.

Consider also an  arbitrary infinite horizon optimal control process associated with the limit of the upper iterates of the Minkowski--Bellman equation. Namely, for any $x$, we construct infinite state and control sequences, $\mathbf{\oz}_\infty(x)=\{\oz_0(x),\oz_1(x),\ldots\}$ and $\mathbf{\ov}_\infty(x)=\{\ov_0(x),\ov_1(x),\ldots\}$, generated, for all $k\in \N$, by
\begin{subequations}
\label{eq:05.08}
\begin{align}
\label{eq:05.08.a}
\oz_{k+1}(x)&\eqbyd A\oz_{k}(x)+B\ov_{k}(x)\text{ with }\oz_{0}(x)\eqbyd x,\\
\label{eq:05.08.b}
\ov_k(x)&\in \ou_\infty(\oz_k(x)),
\end{align}
\end{subequations}
where the selection of controls $\ov_k(x)\in \ou_\infty(\oz_k(x))$ is arbitrary. The limit of the sum of the values of the associated stage cost function $\ell(\oz_j(x),\ov_j(x))$ and the value of the terminal cost function $V_f(\oz_k(x))=\gauge{\mcL}{\oz_k(x)}$ is specified, for all $x\in\R^n$, by 
\begin{align*}
&\oV_\infty(x)=\oV_\infty(\mathbf{\oz}_\infty(x),\mathbf{\ov}_\infty(x))\nonumber\\
&\eqbyd \lim_{k\rightarrow \infty}\left( \sum_{j=0}^{k-1} \ell(\oz_j(x),\ov_j(x))+V_f(\oz_k(x))\right).
\end{align*}
Note that, for any $x\in\R^n$,
\begin{align*}
&\lim_{k\rightarrow \infty}\left( \sum_{j=0}^{k-1} \ell(\oz_j(x),\ov_j(x))+V_f(\oz_k(x))\right)\nonumber\\
&= \lim_{k\rightarrow \infty}\left( \sum_{j=0}^{k-1} \ell(\oz_j(x),\ov_j(x))\right)+\lim_{k\rightarrow \infty}V_f(\oz_k(x)).
\end{align*}
Since, in view of Remark~\ref{rem:04.02}, $\mathbf{\oz}_\infty(x)$ and $\mathbf{\ov}_\infty(x)$ converge to $(0,0)\in \R^n\times\R^m$ for any $x\in \R^n$, we have that $\lim_{k\rightarrow \infty}V_f(\oz_k(x))=\lim_{k\rightarrow \infty}\gauge{\mcL}{\oz_k(x)}=0$ for all $x\in \R^n$, and, in turn, for all $x\in \R^n$,
\begin{align}
\label{eq:05.09}
\oV_\infty(x)= \lim_{k\rightarrow \infty}\left( \sum_{j=0}^{k-1} \ell(\oz_j(x),\ov_j(x))\right).
\end{align}
By construction, $\oV_\infty(x)=\oV_\infty(\mathbf{\oz}_\infty(x),\mathbf{\ov}_\infty(x))$ is guaranteed to converge pointwise in $\R^n$, as $k\rightarrow \infty$ to
\begin{equation}
\label{eq:05.10}
\oV(x)= \gauge{\oP_\infty}{x}.
\end{equation}
Thus, appealing to~\cite[Theorem~1.18.12.]{schneider:1993}, $\oV_\infty\bb$ converges uniformly over the unit sphere $\mathbb{S}^{n-1}$ to $\gaugebb{\oP_\infty}$.

The considered infinite horizon optimal control process $(\mathbf{\oz}_\infty(x),\mathbf{\ov}_\infty(x))$ associated with the limit of the upper iterates of the Minkowski--Bellman equation is not necessarily infinite horizon optimal control process with respect to the limit of the lower iterates of the Minkowski--Bellman equation so that, for all $x\in \R^n$,
\begin{align*}
\uV_\infty(\mathbf{\uz}_\infty(x),\mathbf{\uv}_\infty(x))&\le \uV_\infty(\mathbf{\oz}_\infty(x),\mathbf{\ov}_\infty(x))\nonumber\\
&= \oV_\infty(\mathbf{\oz}_\infty(x),\mathbf{\ov}_\infty(x)).
\end{align*}
Thus, for all $x\in \R^n$, $\gauge{\uP_\infty}{x}\le \gauge{\oP_\infty}{x}$ and, in turn,
\begin{equation}
\label{eq:05.11}
\oP_\infty\subseteq \uP_\infty.
\end{equation}
Likewise, the considered infinite horizon optimal control process $(\mathbf{\uz}_\infty(x),\mathbf{\uv}_\infty(x))$ associated with the limit of the lower iterates of the Minkowski--Bellman equation is  not necessarily infinite horizon optimal control process with respect to the limit of the upper iterates of the Minkowski--Bellman equation so that, for all $x\in \R^n$,
\begin{align*}
\oV_\infty(\mathbf{\oz}_\infty(x),\mathbf{\ov}_\infty(x))&\le \oV_\infty(\mathbf{\uz}_\infty(x),\mathbf{\uv}_\infty(x))\nonumber\\
&= \uV_\infty(\mathbf{\uz}_\infty(x),\mathbf{\uv}_\infty(x)).
\end{align*}
Hence, for all $x\in \R^n$, $\gauge{\oP_\infty}{x}\le \gauge{\uP_\infty}{x}$ and, in turn,
\begin{equation}
\label{eq:05.12}
\uP_\infty\subseteq \oP_\infty.
\end{equation}
Since  $\oP_\infty\subseteq \uP_\infty$ and $\uP_\infty\subseteq \oP_\infty$, the limits of the generator sets of the lower and upper iterates are equal, as formally summarized by the following statement. 
\begin{theorem}
\label{thm:05.01}
Suppose Assumptions~\ref{ass:02.01} and~\ref{ass:02.02} hold. Suppose also that the set $\mcL$ is a proper $C$--set in $\R^n$ that verifies relation~\eqref{eq:04.04}. The set equalities
\begin{equation}
\label{eq:05.13}
\uP_\infty=\oP_\infty
\end{equation}
hold true for the limits $\uP_\infty$ and $\oP_\infty$ of the generator sets associated with the lower and upper iterates of the Minkowski--Bellman equation, respectively.
\end{theorem}
We note that, by Theorem~\ref{thm:05.01}, for all $x\in\R^n$,
\begin{equation}
\label{eq:05.14}
\gauge{\oP_\infty}{x}=\gauge{\uP_\infty}{x}.
\end{equation}

\subsection{Arbitrary Iterates and Independence of Limits}
\label{sec:05.04}

The arbitrary iterates refer to the value functions and optimizer maps generated by the Minkowski--Bellman equation~\eqref{eq:02.06} with the boundary condition~\eqref{eq:02.07}. As in the previously considered cases, the solution of  the Minkowski--Bellman equation~\eqref{eq:02.06} with the boundary condition~\eqref{eq:02.07} is entirely determined by the associated generator sets $\mcP_k$ produced by the set--dynamics~\eqref{eq:04.01} with the initial condition in its polar form~\eqref{eq:04.02}. 

The structural and boundedness properties of the generator sets $\mcP_k$ and the related value functions, optimizer maps and selections, $V_k\bb$, $u_k\bb$ and $\nu_k\bb$, are identical to those discussed in the previously considered cases. However, Proposition~\ref{prop:04.03} and Theorem~\ref{thm:04.01} do not apply to this setting. Hence, the monotonicity, convergence and fixed point properties  need additional discussion. 

Under Assumptions~\ref{ass:02.01},~\ref{ass:02.02} and~\ref{ass:02.03}, the existence of a proper $C$--set $\mcL$ in $\R^n$ that satisfies both relations~\eqref{eq:04.04} and~\eqref{eq:04.08} is guaranteed by Proposition~\ref{prop:04.01}. Let the initial condition of the upper iterates be such a set $\mcL$. In this setting, Proposition~\ref{prop:04.04} yields, for all $k\in \N$,
\begin{equation}
\label{eq:05.15}
\oP_{k+1}\subseteq \mcP_{k+1}\subseteq \uP_{k+1},
\end{equation}
where the sets $\uP_k$ and $\oP_k$ denote the generator sets of the corresponding lower and upper iterates of the Minkowski--Bellman equation. The relations~\eqref{eq:05.15} imply that, in terms of the related value functions, for all $k\in \N$ and all $x\in \R^n$, 
\begin{equation}
\label{eq:05.16}
\uV_{k+1}(x)\le V_{k+1}(x)\le \oV_{k+1}(x),
\end{equation}
where value functions $\uV_k\bb$ and $\oV_k\bb$ denote the value functions of the corresponding lower and upper iterates of the Minkowski--Bellman equation.

Proposition~\ref{prop:04.04} and Theorems~\ref{thm:04.01},~\ref{thm:04.02} and~\ref{thm:05.01} guarantee that the sequence of the generator sets $\mcP_k$ converges to a proper $C$--set $\mcP_\infty$ in $\R^n$. In particular, Theorem~\ref{thm:04.01} yields convergence of the sequences of the generator sets $\uP_k$ as well as the generator sets $\oP_k$. Theorem~\ref{thm:05.01} has established the set equality $\oP_\infty=\uP_\infty$ of the respective limits. With these facts and relations~\eqref{eq:05.15} established in Proposition~\ref{prop:04.04} in mind, the convergence of the sequences of the generator sets $\mcP_k$ to proper $C$--set $\mcP_\infty$ in $\R^n$ is guaranteed. Indeed, Theorem~\ref{thm:04.02} guarantees the relations
\begin{equation}
\label{eq:05.17}
\oP_\infty\subseteq \mcP_\infty\subseteq \uP_\infty
\end{equation}
for the limits of the generator sets of the lower, arbitrary and upper iterates of the Minkowski--Bellman equation. Theorem~\ref{thm:05.01} ensures $\oP_\infty=\uP_\infty$. Hence, the limits of the generator sets of the lower, arbitrary and upper iterates of the Minkowski--Bellman equation are identical. Theorems~\ref{thm:04.01} and~\ref{thm:05.01} yield additional fixed point properties of this limit. These facts are summarized by the following.
\begin{theorem}
\label{thm:05.02}
Suppose Assumptions~\ref{ass:02.01},~\ref{ass:02.02} and~\ref{ass:02.03} hold. Suppose also that the set $\mcL$ is a proper $C$--set in $\R^n$ that verifies relations~\eqref{eq:04.04} and~\eqref{eq:04.08}.
\begin{itemize}
\item [(i)] The sequence $\{\mcP_k\}_{k\ge 1}$ of the generator sets $\mcP_k$ converges to a proper $C$--set $\mcP_\infty$ in $\R^n$.
\item [(ii)] The set equalities 
\begin{equation}
\label{eq:05.18}
\oP_\infty=\mcP_\infty=\uP_\infty
\end{equation}
hold true for the limits $\uP_\infty$,  $\mcP_\infty$ and $\oP_\infty$ of the generator sets of the corresponding lower, arbitrary and upper iterates of the Minkowski--Bellman equation.
\item [(iii)] The sequence $\{\mcT_k\}_{k\ge 1}$ of the sets $\mcT_k$ generated by~\eqref{eq:04.09} converges to a proper $C$--set $\mcT_\infty$ in $\R^{n+m}$, which satisfies relation~\eqref{eq:04.20}.
\item [(iv)] The limit $\mcP_\infty$ solves the fixed point set--equation~\eqref{eq:04.03}.
\end{itemize}  
\end{theorem}

\subsection{Uniqueness and Stability of Fixed Point}
\label{sec:05.05}

Theorem~\ref{thm:05.02} guarantees that the fixed point set--equation~\eqref{eq:04.03} admits a unique solution over the space of proper $C$--sets in $\R^n$. The boundedness and convergence properties, guaranteed due to Propositions~\ref{prop:04.01}--\ref{prop:04.05} and Theorems~\ref{thm:04.01}--\ref{thm:05.02}, imply directly that this unique fixed point is an asymptotically stable attractor for the set--dynamics~\eqref{eq:03.07}. The domain of attraction, specified in terms of related polar sets, is the space of $C$--sets in $\R^n$. 
\begin{theorem}
\label{thm:05.03}
Suppose Assumptions~\ref{ass:02.01} and~\ref{ass:02.02} hold.   There exists a unique proper $C$--set $\mcP_\infty$ in $\R^n$ solving the fixed point set--equation~\eqref{eq:04.03} 
\begin{align*}
\mcP=P_x\left(\mcC^*\oplus (A\ B)^T\mcP^*\right)^*.
\end{align*}
Furthermore, the sets  $\mcP_\infty$ is an asymptotically stable attractor for the set--dynamics~\eqref{eq:03.07} 
\begin{align*}
\mcP^+=P_x\left(\mcC^*\oplus (A\ B)^T\mcP^*\right)^*.
\end{align*}
The domain of attraction, in terms of polar sets $\mcP^*$, is the space of $C$--sets in $\R^n$. 
\end{theorem}
Since Theorem~\ref{thm:05.03} applies when $\mcP_0^*=\{0\}$, the fixed point $\mcP_\infty$ is unique over the space of nonempty convex closed subsets of $\R^n$ containing the origin as an interior point. 

Theorem~\ref{thm:05.03} guarantees that the fixed point of the Minkowski--Bellman equation~\eqref{eq:02.08} is unique, in terms of the value function, over the space of Minkowski functions of proper $C$--sets in $\R^n$. This unique fixed point is formed by the value function $\gaugebb{\mcP_\infty}$, while its generator set $\mcP_\infty$ and associated set $\mcT_\infty=\left(\mcC^*\oplus (A\ B)^T\mcP_\infty^*\right)^*$ entirely characterize the fixed point optimizer map $u_\infty\bb$ and related selection $\nu_\infty\bb$, which are, respectively, a $\Pi$--class set--valued map and a $\pi$--class single--valued function. Namely,  $\gaugebb{\mcP_\infty}$ is the unique Minkowski function of proper $C$--sets in $\R^n$ such that, for  all $x\in \R^n$ and all $u\in u_\infty(x)$ (including $u=\nu_\infty(x)$),
\begin{align}
\gauge{\mcP_\infty}{x}&=\gauge{\mcC}{(x,u)}+\gauge{\mcP_\infty}{Ax+Bu},\text{ where}\nonumber\\
u_\infty(x)&=\arg\min_{u}\ \gauge{\mcC}{(x,u)}+\gauge{\mcP_\infty}{Ax+Bu}\nonumber\\
&=\{u\in \R^m\ :\ (x,u)\in \gauge{\mcP_\infty}{x}\mcT_\infty\}\text{ and}\nonumber\\
\label{eq:05.19}
\nu_\infty(x)&=\arg\min_u\{u^Tu\ :\ u\in u_\infty(x)\}.
\end{align}
Any arbitrary iteration of the Minkowski--Bellman equation~\eqref{eq:02.06} with the boundary condition~\eqref{eq:02.07} converges asymptotically  in a stable manner to the established fixed point. The convergence of the value functions, optimizer maps and selections, $V_k\bb$, $u_k\bb$ and $\nu_k\bb$, is uniform over the unit sphere $\mathbb{S}^{n-1}$, which implies both pointwise convergence over $\R^n$ and uniform convergence over nonempty compact subsets of $\R^n$. The fixed point value function, optimizer map and selection are $\gaugebb{\mcP_\infty}$, $u_\infty\bb$ and $\nu_\infty\bb$.

\section{Polytopic Iterates}
\label{sec:06}

\subsection{Refined Properties}
\label{sec:06.01}

Here, we elaborate on the case when the generator sets $\mcC$ and $\mcQ_f$ of the stage and terminal cost functions $\ell\bbb$ and $V_f\bb$ are proper $C$--polytopes. In the proper $C$--polytopic setting, it is worth observing that: $(i)$ the polar set of a proper $C$--polytope set is itself a proper $C$--polytope, $(ii)$ Minkowski sum of a proper $C$--polytope and a $C$--polytope is a proper $C$--polytope, $(iii)$ linear transformation from $\R^n$ to $\R^{n+m}$ of proper $C$--polytope in $\R^n$ is a $C$--polytope in $\R^{n+m}$, and $(iv)$ the $(x,u)\mapsto x$ projection of a proper $C$--polytope in $\R^{n+m}$ is a proper $C$--polytope in $\R^n$. Hence, in this setting, the iterates of the set--dynamics~\eqref{eq:04.01} preserve the proper $C$--polytopic structure of the initial condition $\mcP_0:=\mcQ_f$. Consequently, the generator sets $\mcP_k$ of the iterates of the associated Minkowski--Bellman equation as well as related sets $\mcT_k$ specified in~\eqref{eq:04.09} are proper $C$--polytopes for all $k\in \N$. The limiting generator set $\mcP_\infty$ and related limiting set $\mcT_\infty$ can not be guaranteed to be proper $C$--polytopes; Rather these limiting sets can be guaranteed to be proper $C$--sets. However, in the case of finite determination, i.e. when $\mcP_{k^*+1}=\mcP_{k^*+2}$ holds for a finite integer $k^*$, the limiting sets $\mcP_\infty=\mcP_{k^*+2}$ and $\mcT_\infty=\mcT_{k^*+2}$ are guaranteed to be proper $C$--polytopes. 

This setting allows for refinement of the structural properties of the value functions $V_k\bb$, optimizer maps $u_k\bb$ and selections $\nu_k\bb$. The value functions $V_k\bb$ are Minkowski functions of proper $C$--polytopes $\mcP_k$ in $\R^n$. Consequently, the value functions are additionally piecewise linear. More specifically, the generator sets $\mcP_k$ admit an irreducible representation given by
\begin{equation}
\label{eq:06.01}
\mcP_k=\{x\ :\ \forall i\in\mcI_k,\ \alpha_{(k,i)}^Tx\le 1\},
\end{equation}
where $\mcI_k$ is a finite index set, while the collection of vectors $\{ \alpha_{(k,i)}\in \R^n\ :\ i\in \mcI_k\}$ spans $\R^n$ and it also induces a conical partition $\{\mcR_{(k,i)}\ :\ i\in\mcI_k\}$ of $\R^n$, in which each of the cones $\mcR_{(k,i)},\ i\in\mcI_k$ is given by
\begin{equation}
\label{eq:06.02}
\mcR_{(k,i)}:=\{x\ :\ \forall j\in\mcI_k,\ (\alpha_{(k,j)}-\alpha_{(k,i)})^Tx\le 0\}.
\end{equation}
The related value functions $V_k\bb$ satisfy, for all $x\in \R^n$, 
\begin{align}
V_k(x)&= \max_{i\in\mcI_k}\alpha_{(k,i)}^Tx\text{ so that}\nonumber\\
\label{eq:06.03}
V_k(x)&= \alpha_{(k,i)}^Tx\text{ when }x\in \mcR_{(k,i)}.
\end{align}
Likewise, the related sets $\mcT_k$ admit an irreducible representation given by
\begin{equation}
\label{eq:06.04}
\mcT_k=\{(x,u)\ :\ \forall j\in\mcJ_k,\ \beta_{(k,j)}^Tx+\gamma_{(k,j)}^Tu\le 1\},
\end{equation}
where $\mcJ_k$ is a finite index set and the collection of vectors $\{ (\beta_{(k,j)},\gamma_{(k,j)})\in \R^{n+m}\ :\ j\in \mcJ_k\}$ spans $\R^{n+m}$.
The associated optimizer maps $u_k\bb$ satisfy, for all $x\in \R^n$, 
\begin{align}
u_k(x)=\{u\ :\ &\forall j\in\mcJ_k,\ \gamma_{(k,j)}^Tu\le (\alpha_{(k,i)}-\beta_{(k,j)})^Tx\}\nonumber\\
\label{eq:06.05}
&\text{when }x\in \mcR_{(k,i)}.
\end{align}
The related selections $\nu_k\bb$ are given, for all $x\in \R^n$, by
\begin{align}
\nu_k(x)=\arg\min_u\{u^Tu\ :\ &\forall j\in\mcJ_k,\nonumber\\
 &\gamma_{(k,j)}^Tu\le (\alpha_{(k,i)}-\beta_{(k,j)})^Tx\}\nonumber\\
\label{eq:06.06}
&\text{when }x\in \mcR_{(k,i)}.
\end{align}
Hence, the optimizer maps $u_k\bb$ are polyhedral set--valued maps (i.e. their graphs are unions of finitely many polyehdral sets) and are actually polytopic--valued (for all $x\in\R^n$, the set $u_k(x)$ is a polytope). The optimizer maps $u_k\bb$ are positively homogeneous of the first degree. The polytopic structure also allows for a refinement of the continuity properties of the optimizer maps and related selections. In view of results related to the Lipschitz continuity of polyhedral set--valued maps~\cite{walkup:wets:1969:a,walkup:wets:1969:b,robinson:1981,klatte:gisbert:1995}, the optimizer maps  $u_k\bb$ are Lipschitz continuous with respect to the Hausdorff distance. When the optimizer maps $u_k\bb$ are single--valued, they are piecewise linear $\pi$--class functions, as implied by~\eqref{eq:06.05}. Furthermore, the associated selections $\nu_k\bb$ are piecewise linear $\pi$--class single--valued functions, as dictated by~\eqref{eq:06.06}. 

The fixed point value function $V_\infty\bb=\gaugebb{\mcP_\infty}$ is guaranteed to be the Minkowski function of a proper $C$--set $\mcP_\infty$ in $\R^n$. In the case of finite determination, the fixed point value function $V_\infty\bb=\gaugebb{\mcP_\infty}=\gaugebb{\mcP_{k^*+2}}$ is the Minkowski function of a proper $C$--polytope $\mcP_\infty=\mcP_{k^*+2}$ in $\R^n$. The fixed point optimizer map $u_\infty\bb$ is guaranteed to be a $\Pi$--class set--valued map, which is a $\pi$--class function when it is single--valued. In the case of finite determination, the fixed point optimizer map $u_\infty\bb=u_{k^*+2}\bb$ is a polytopic--valued, positively homogeneous of the first degree, polyhedral and Lipschitz continuous, with respect to the Hausdorff distance, set--valued map, which is a  piecewise linear $\pi$--class function when it is single--valued. 
The fixed point selection $\nu_\infty\bb$ is guaranteed to be a $\pi$--class single--valued function. In the case of finite determination, the fixed point selection $\nu_\infty\bb=\nu_{k^*+2}\bb$ is, in fact, a  piecewise linear $\pi$--class single--valued function. 

\subsection{Approximations}
\label{sec:06.02}

The polytopic computations can be used to construct lower and upper proper $C$--polytopic approximations of the generator sets of arbitrary iterates of the Minkowski--Bellman equation and its fixed point. Namely, \cite[Theorem~1.8.13.]{schneider:1993} guarantees that for any given pair of proper $C$--sets $\mcC$ and $\mcQ_f$ in $\R^{n+m}$ and $\R^n$ and any $\delta> 0$, we can construct pairs of proper $C$--polytopes, $\uC$ and $\oC$ and $\uQ_f$ and $\oQ_f$, in $\R^{n+m}$ and $\R^n$ such that $H_{\mcB^{n+m}}(\uC,\oC)\le \delta$, $H_{\mcB^n}(\uQ_f,\oQ_f)\le \delta$ and
\begin{equation*}
\oC\subseteq \mcC\subseteq \uC\text{ and }\oQ_f\subseteq \mcQ_f\subseteq \uQ_f
\end{equation*}
Connection of these relations with the solution of the Minkowski--Bellman equation, for which $\mcC$ and $\mcQ_f$ are the generator sets of the stage and terminal cost functions $\ell\bbb$ and $V_f\bb$ is clear. Namely, the solution of the Minkowski--Bellman equation, for which $\uC$ and $\uQ_f$ are employed as the generator sets of the stage and terminal cost functions $\ell\bbb$ and $V_f\bb$ yields polytopic lower approximation of the solution to the original Minkowski--Bellman equation. Likewise, the solution of the Minkowski--Bellman equation, for which $\oC$ and $\oQ_f$ are utilized as the generator sets of the stage and terminal cost functions $\ell\bbb$ and $V_f\bb$ yields a polytopic upper approximation of the solution to the original Minkowski--Bellman equation. Thus, such polytopic iterates can be used to lower and upper approximate the arbitrary iterates of the Minkowski--Bellman equation as well as its fixed point. The lower and upper polytopic approximations of the fixed point of the Minkowski--Bellman equation are provided by the lower and upper polytopic iterates obtained for large enough $k$. The quality of such lower and upper approximations can be regulated via selection of  $\delta>0$. 

A further computational convenience can be ensured by utilizing~\cite[Theorem~1.8.15.]{schneider:1993}, in view of which, for any given pair of proper $C$--sets $\mcC$ and $\mcQ_f$ in $\R^{n+m}$ and $\R^n$ and any $\delta> 0$, we can construct a pair of proper $C$--polytopes, $\widetilde{\mcC}$ and $\widetilde{\mcQ}_f$ in $\R^{n+m}$ and $\R^n$ such that
\begin{equation*}
\widetilde{\mcC}\subseteq \mcC\subseteq (1+\delta)\widetilde{\mcC}\text{ and }\widetilde{\mcQ}_f\subseteq \mcQ_f\subseteq (1+\delta)\widetilde{\mcQ}_f.
\end{equation*}
Within this construction, the solution of the Minkowski--Bellman equation, for which $\widetilde{\mcC}$ and $\widetilde{\mcQ}_f$ are employed as the generator sets of the stage and terminal cost functions $\ell\bbb$ and $V_f\bb$ can be used to construct both lower and upper polytopic approximations of the solution to the original Minkowski--Bellman equation, its iteration and its fixed point. The quality of the related approximations can be controlled via choice of $\delta >0$.

\section{Conclusions}
\label{sec:07}

It has been established that the Minkowski--Bellman equation, its iteration and its fixed point  are well posed. The characterization of the corresponding value functions, optimizer maps and their selections has been derived. In particular, it has been demonstrated that the related value functions are Minkowski function of clearly defined proper $C$--sets  as well as that the related optimizer maps and their selections are, respectively, $\Pi$--class set--valued maps and $\pi$--class single--valued functions. These properties have been further refined in the proper $C$--polytopic setting.

\subsection*{Acknowledgements} 

The author is grateful to Zvi Artstein for helpful feedback on the continuity aspects for optimizer maps $u_k\bb$.

\subsection*{APPENDIX A: Compact Form of Stage Cost}
When, for all $(x,u)\in \R^n\times \R^m$,
\begin{equation*}
\ell(x,u)=\gauge{\mcQ}{x}+\gauge{\mcS}{(x,u)}+\gauge{\mcR}{u},
\end{equation*}
for proper $C$--sets $\mcQ$, $\mcS$ and $\mcR$, the proper $C$--set $\mcC$ in the compact form of the stage cost function $\ell\bbb$ of~\eqref{eq:02.02} is
\begin{equation*}
\mcC:=\left(P_x^T\mcQ^*\oplus \mcS^*\oplus P_u^T\mcR^*\right)^*.
\end{equation*}
If $\gauge{\mcS}{(x,u)}$ is absent, i.e., for all $(x,u)\in \R^n\times \R^m$, $\ell(x,u)=\gauge{\mcQ}{x}+\gauge{\mcR}{u}$ then $\mcC:=\left(P_x^T\mcQ^*\oplus P_u^T\mcR^*\right)^*$.

\subsection*{APPENDIX B--1: Proof of Theorem~\ref{thm:03.04}}
$(i)$: $\mcC^*$ is a proper $C$--set in $\R^{n+m}$ since $\mcC$ is a proper $C$--set in $\R^{n+m}$. By the same token, $\mcP^*$ is a proper $C$--set in $\R^n$. Thus, $(A\ B)^T\mcP^*$ is a $C$--set in $\R^{n+m}$. Hence, $\mcC^*\oplus (A\ B)^T\mcP^*$ is a proper $C$--set in $\R^{n+m}$,  and, in turn, $\mcT^+= \left(\mcC^*\oplus (A\ B)^T\mcP^*\right)^*$ is a proper $C$--set in $\R^{n+m}$. $\mcP^+=P_x\mcT^+$ is a proper $C$--set in $\R^n$ being a $(x,u)\mapsto x$ projection of a proper $C$--set $\mcT^+$ in $\R^{n+m}$ to $\R^n$.

Before proceeding, we note that a direct calculation yields, for all $(x,u)\in\R^n\times\R^m$,
\begin{align*}
J^+(x,u)&=\ell(x,u)+V(Ax+Bu)\\
&=\gauge{\mcC}{(x,u)}+\gauge{\mcP}{Ax+Bu}\\
&=\support{\mcC^*}{(x,u)}+\support{\mcP^*}{Ax+Bu}\\
&=\support{\mcC^*}{(x,u)}+\support{\mcP^*}{(A\ B)(x,u)}\\
&=\support{\mcC^*}{(x,u)}+\support{(A\ B)^T\mcP^*}{(x,u)}\\
&=\support{ \mcC^*\oplus (A\ B)^T\mcP^*}{(x,u)}\\
&=\support{(\mcT^+)^*}{(x,u)}\\
&=\gauge{\mcT^+}{(x,u)},
\end{align*}
so that, for all $x\in \R^n$,
\begin{align*}
V^+(x)&=\min_u\ J^+(x,u),\text{ and}\\
u^+(x)&=\arg\min_u\ J^+(x,u).
\end{align*}
Under invoked assumptions, the values of the objective and value functions, $J^+(x,u)$ and $V^+(x)$, are nonnegative and finite for all $(x,u)\in\R^n\times \R^m$ and all $x\in\R^n$, while the related $\arg\min$ set is nonempty for all $x\in \R^n$. 

$(ii)$: $V^+\bb$ is a sublinear function:\\
For any $x_1\in \R^n$ and $x_2\in\R^n$  let $u_1$ and $u_2$ denote any of related optimizer points, so that $V^+(x_1)=J^+(x_1,u_1)$ and $V^+(x_2)=J^+(x_2,u_2)$. Hence, since $J^+\bbb$ is a sublinear function, 
\begin{align*}
\min_uJ^+(x_1+x_2,u)&\le J^+(x_1+x_2,u_1+u_2)\\
&\le J^+(x_1,u_1)+J^+(x_2,u_2),\text{ i.e.}\\
V^+(x_1+x_2)&\le V^+(x_1)+V^+(x_2)
\end{align*}
and $V^+\bb$ is a subadditive function. For any $x\in \R^n$ and $\eta=0$, $V^+(\eta x)=V^+(0 x)=V^+(0)=J^+(0,0)=0$ and $\eta V^+(x)=0 V^+(x)=0$. So, for all $x\in \R^n$ and $\eta=0$, $V^+(\eta x)=\eta V^+(x)$. For any $x\in \R^n$, let $u^0$ be any of related optimizer points, so that $V^+(x)=J^+(x,u^0)$. Then, since $J^+\bbb$ is a sublinear function, for any $\eta >0$, 
\begin{align*}
\min_uJ^+(\eta x,u)&\le J^+(\eta x,\eta u^0)=\eta J^+(x,u^0),\text{ i.e.}\\
V^+(\eta x)&\le\eta V^+(x).
\end{align*}
Clearly, $V^+(\eta x)=J^+(\eta x,u_\eta)=\eta J^+(x,\eta^{-1}u_\eta)$ for some optimizer point $u_\eta$ at $x_\eta:=\eta x$. If $V^+(\eta x)< \eta V^+(x)$, then  $\eta J^+(x,\eta^{-1}u_\eta)<\eta V^+(x)$ or, equivalently, $J^+(x,\eta^{-1}u_\eta)< V^+(x)$. This contradicts optimality of $V^+(x)$. Thus, for all $x\in \R^n$ and all $\eta\ge 0$, 
\begin{equation*}
V^+(\eta x)=\eta V^+(x).
\end{equation*}
Hence,  $V^+\bb$ is a sublinear function.

$\forall x\in\R^n,\ V^+(x)=\gauge{\mcP^+}{x}$:\\
Since $\mcT^+$ is a proper $C$--set in $\R^{n+m}$, $J^+(x,u)>0$ for all $(x,u)\in \R^n\times \R^m\setminus\{(0,0)\}$ and $V^+(x)>0$ for all $x\in\R^n\setminus\{0\}$. Hence, due to Theorem~\ref{thm:03.02}, there is a unique proper $C$--set in $\R^n$ with support function $V^+\bb$. In turn, by Theorem~\ref{thm:03.03}, there is a unique proper $C$--set in $\R^n$ whose Minkowski function is $V^+\bb$. Indeed, $V^+\bb$ is the Minkowski function of the proper $C$--set $\mcP^+$ in $\R^n$.

Take arbitrary $x\in\R^n$ and let $\rho=V^+(x)$. Then $V^+(x)=J^+(x,u^0)=\gauge{\mcT^+}{(x,u^0)}$ for some $u^0\in \R^m$. In turn, $(x,u^0)\in \rho \mcT^+$. Consequently, $x\in \rho P_x\mcT^+=\rho \mcP^+$ and $\gauge{\mcP^+}{x}\le \rho$. If $\gauge{\mcP^+}{x}=\theta< \rho$ then there exists a $v\in\R^m$ such that $(x,v)\in \theta \mcT^+$ and, in turn, $J^+(x,v)=\gauge{\mcT^+}{(x,v)}\le \theta <\rho=V^+(x)$. This contradicts the optimality of $V^+\bb$. Hence, $\gauge{\mcP^+}{x}=\rho$.

Take arbitrary $x\in \R^n$ and let $\rho=\gauge{\mcP^+}{x}$. There is a $v\in \R^m$ such that $(x,v)\in \rho \mcT^+$ and $V^+(x)\le J^+(x,v)=\gauge{\mcT^+}{(x,v)} \le \rho$. If $V^+(x)=\theta <\rho$ then there exists a $u^0\in \R^m$ such that $V^+(x)=J^+(x,u^0)=\gauge{\mcT^+}{(x,u^0)}$. Thus, $(x,u^0)\in \theta \mcT^+$ so that $x\in \theta \mcP^+$ and $\gauge{\mcP^+}{x}\le \theta <\rho$ contradicting $\rho=\gauge{\mcP^+}{x}$. Hence, $V^+(x)=\rho$.

Summa summarum, for all $x\in\ \R^n,\ V^+(x)=\gauge{\mcP^+}{x}$.

$(iii)$: Since $\mcT^+$ is a proper $C$--set in $\R^n$, $J^+\bbb=\gauge{\mcT^+}{(\cdot,\cdot)}$ satisfies postulates of~\cite[Theorem~1.17., Proposition~2.22., and Theorem~7.41.]{rockafellar:wets:2009}. Hence, $u^+(x)=\arg\min_u\ J^+(x,u)$ is: $(a)$ nonempty and compact for all $x\in \R^n$ by~\cite[Theorem~1.17.]{rockafellar:wets:2009}, $(b)$  convex for all $x\in\R^n$ by~\cite[Proposition~2.22.]{rockafellar:wets:2009}, and $(c)$ locally bounded and outer semicontinuous for all $\R^n$ by~\cite[Theorem~7.41.]{rockafellar:wets:2009}.  

$\forall x\in\R^n,\ u^+(x)=\{u\in\R^m\ :\ (x,u)\in \gauge{\mcP^+}{x}\mcT^+\}$:\\
Take any $x\in \R^n$ and any $v\in \{u\in\R^m\ :\ (x,u)\in \gauge{\mcP^+}{x}\mcT^+\}$. Then $(x,v)\in \gauge{\mcP^+}{x}\mcT^+$ and $J^+(x,v)= \gauge{\mcT^+}{(x,v)}\le \gauge{\mcP^+}{x}=V^+(x)$. If $J^+(x,v)< V^+(x)$ optimality of $V^+(x)$ is contradicted. So, it must hold $V^+(x)=J^+(x,v)$, and $v\in u^+(x)$. Hence, $\{u\in\R^m\ :\ (x,u)\in \gauge{\mcP^+}{x}\mcT^+\}\subseteq u^+(x)$.

Take any $x\in \R^n$ and any $u^0\in u^+(x)$. Then $\gauge{\mcP^+}{x}=V^+(x)=J^+(x,u^0)=\gauge{\mcT^+}{(x,u^0)}$. In turn, $(x,u^0)\in \gauge{\mcP^+}{x}\mcT^+$ and $u^0\in \{u\in\R^m\ :\ (x,u)\in \gauge{\mcP^+}{x}\mcT^+\}$. Hence, $u^+(x)\subseteq \{u\in\R^m\ :\ (x,u)\in \gauge{\mcP^+}{x}\mcT^+\}$.

Thus, for all $x\in\R^n,\ u^+(x)=\arg\min_u\ J^+(x,u)=\{u\in\R^m\ :\ (x,u)\in \gauge{\mcP^+}{x}\mcT^+\}$.

$u^+\bb$ is positively homogeneous of the first degree:\\
Take any $x\in \R^n$ and any $\eta\in \R_{\ge 0}$. For $\eta=0$, $\eta x=0$ and $u^+(\eta x)=u^+(0)=\{0\}$.  Also, $u^+(x)$ is compact so that $\eta u^+(x)=0u^+(x)=\{0\}$. So that $u^+(\eta x)=\eta u^+(x)$ for $\eta =0$. Note that we used the fact that, by definition, $u^+(0)=\{0\}$. Consider from now on $\eta >0$. Take any $v\in u^+(\eta x)$. Then, $(\eta x,v)\in \gauge{\mcP^+}{\eta x} \mcT^+=\eta \gauge{\mcP^+}{x} \mcT^+$. Since $\mcT^+$ is a proper $C$--set in $\R^{n+m}$ and $\eta >0$, we have
\begin{align*}
\eta^{-1}(\eta x,v)=(x,\eta^{-1}v)&\in \gauge{\mcP^+}{x} \mcT^+,\text{ or, equivalently,}\\
\eta^{-1}v&\in u^+(x),\text{ i.e. }v \in \eta u^+(x).
\end{align*}
Thus, $u^+(\eta x)\subseteq \eta u^+(x)$.\\
Next, take any $v\in \eta u^+(x)$. Then, $\eta^{-1}v\in u^+(x)$ and, in turn, $(x,\eta ^{-1}v)\in \gauge{\mcP^+}{x}\mcT^+$. Since $\mcT^+$ is a proper $C$--set in $\R^{n+m}$ and $\eta >0$, it follows that 
\begin{equation*}
\eta(x, \eta^{-1}v)=(\eta x,v)\in \eta\gauge{\mcP^+}{x} \mcT^+=\gauge{\mcP^+}{\eta x} \mcT^+,
\end{equation*}
or, equivalently, 
\begin{equation*}
v\in u^+(\eta x).
\end{equation*}
Thus, $\eta u^+(x)\subseteq u^+(\eta x)$.  Hence, $u^+(\eta x)=\eta u^+(x)$ for all $x\in \R^n$ and all $\eta\in\R_{\ge 0}$, and $u^+\bb$ is positively homogeneous of the first degree.

Summa summarum, $u^+\bb\ :\ \R^n\rightrightarrows \R^m$ is a $\Pi$--class set--valued map such that, for all $x\in \R^n$, 
\begin{equation*}
u^+(x)=\{u\in\R^m\ :\ (x,u)\in \gauge{\mcP^+}{x}\mcT^+\}.
\end{equation*}
$(iv)$: $\nu\bb$ is single--valued for all $x\in\R^n$:\\
For all $x\in\R^n$, $u\mapsto u^Tu$ is strictly convex in $u$ and $u^+(x)$ is a convex and compact subset of $\R^m$. Thus, $\nu^+(x)=\arg\min_{u}\{u^Tu\ :\ u\in u^+(x)\}$ is singleton for all $x\in\R^n$.

$\nu\bb$ is positively homogeneous of the first degree:\\
For all $x\in \R^n$ with $u=\nu^+(x)$, and all $\eta \ge 0$,
\begin{align*}
&(\eta u)^T(\eta u)\le (\eta v)^T(\eta v)\text{ for all }v\in u^+(x),\text{ i.e.}\\ 
&(\eta u)^T(\eta u)\le w^Tw\text{ for all }w\in \eta u^+(x)=u^+(\eta x).  
\end{align*} 
Hence, for all $x\in \R^n$ and all $\eta \ge 0$, $\nu^+(\eta x)=\eta \nu^+(x)$, and $\nu^+\bb$ is positively homogeneous of the first degree.

$\nu\bb$ is continuous:\\
Since $\gaugebb{\mcP^+}$ is the Minkowski function of a proper $C$--set $\mcP^+$ in $\R^n$ and $\mcT^+$ is a proper $C$--set in $\R^n$, $\nu\bb$ is locally bounded and outer semicontinuous for all $x\in \R^n$ by virtue of~\cite[Theorem~7.42.]{rockafellar:wets:2009}. But, $\nu\bb$ is single--valued for all $x\in\R^n$. Hence, $\nu\bb$ is, in fact, continuous for all $x\in\R^n$ by~\cite[Corollary~5.20.]{rockafellar:wets:2009}. 

Summa summarum, $\nu^+\bb\ :\ \R^n\rightarrow \R^m$ is a $\pi$--class single--valued function that is, by its definition, such that, for all $x\in \R^n$, 
\begin{equation*}
\nu^+(x)\in u^+(x).
\end{equation*}
\subsection*{APPENDIX B--2: Proof of Proposition~\ref{prop:04.01}}
$(i)$: Let $u=Kx$ where $K\in\Rnm{m}{n}$ is such that $(A+BK)$ is strictly stable (i.e. the spectral radius of $(A+BK)$ is strictly less than $1$). Note that the function $x\mapsto\ell(x,Kx)$ is sublinear so that, in light of the fact that $\mcC$ is a proper $C$--set in $\R^{n+m}$ and Theorem~\ref{thm:03.02}, it is the Minkowski function of a unique proper $C$--set denoted by $\mcC_x$ in $\R^n$, i.e., for all $x\in \R^n$,
\begin{equation*}
\gauge{\mcC_x}{x}=\gauge{\mcC}{(x,Kx)}.
\end{equation*}
By~\cite[Theorem~1.]{rakovic:2017}, the set $\mcL$ given by
\begin{equation*}
\mcL\eqbyd \left(\bigoplus_{k=0}^\infty ((A+BK)^T)^k\mcC_x^*\right)^*
\end{equation*}
is a proper $C$--set in $\R^n$ verifying the desired relation with equality under $u=Kx$, i.e. such that, for all $x\in \R^n$,
\begin{equation*}
\gauge{\mcL}{(A+BK)x}+\gauge{\mcC}{(x,Kx)}=\gauge{\mcL}{x}.
\end{equation*}
$(ii)$: In view of Theorem~\ref{thm:03.04} and its proof, $\mcK^+$ is a proper $C$--set in $\R^{n+m}$, and, for all $(x,u)\in \R^n\times \R^m$, 
\begin{equation*}
\gauge{\mcK^+}{(x,u)}=\gauge{\mcL}{Ax+Bu}+\gauge{\mcC}{(x,u)}.
\end{equation*}
By $(i)$, for all $x\in\R^n$ there exists a $u\in \R^m$ such that
\begin{equation*}
\gauge{\mcK^+}{(x,u)}\le \gauge{\mcL}{x}.
\end{equation*}
In light of Theorem~\ref{thm:03.04} and its proof, for all $x\in \R^n$,
\begin{equation*}
\gauge{\mcL^+}{x}=\min_u \gauge{\mcK^+}{(x,u)}, 
\end{equation*}
and, thus,  for all $x\in\R^n$ there exists a $u\in \R^m$ such that
\begin{equation*}
\gauge{\mcL^+}{x}\le \gauge{\mcK^+}{(x,u)}\le \gauge{\mcL}{x}. 
\end{equation*}
$(iii)$: By definition, $\mu^*\eqbyd\min_\eta\{\eta\ :\ \mcL\subseteq \eta\mcQ_f,\ \eta\ge 0\}$ 
is strictly positive and finite. Furthermore, for all $\mu\ge \mu^*$, $\mu^{-1}\mcL\subseteq \mcQ_f\text{ or, equivalently, }\mcL\subseteq \mu\mcQ_f$.
Take any proper $C$--set $\mcL$ in $\R^n$ verifying relation~\eqref{eq:04.04} in $(i)$, so that, for all $x\in \R^n$, there exists a $u\in \R^m$ such that, for all $\mu\ge 1$, 
\begin{equation*}
\mu\gauge{\mcL}{Ax+Bu}+\gauge{\mcC}{(x,u)}\le \mu\gauge{\mcL}{x}.
\end{equation*}
In turn, for all $x\in \R^n$, there exists a $u\in \R^m$ such that, for all $\mu\ge \max\{\mu^*,1\}$ 
\begin{align*}
\gauge{\mu^{-1}\mcL}{Ax+Bu}+\gauge{\mcC}{(x,u)}\le \gauge{\mu^{-1}\mcL}{x},
\end{align*}
with
\begin{equation*}
\mu^{-1}\mcL\subseteq \mcQ_f.
\end{equation*}
Thus, any of proper $C$--sets in $\R^n$ of the form $\mu^{-1}\mcL$ with $\mu\in [\max\{\mu^*,1\},\infty)$ verifies the claim.

\subsection*{APPENDIX B--3: Proof of Proposition~\ref{prop:04.02}}
If for some $j\in \N$, $\mcP_{j+1}$ is a proper $C$--set in $\R^n$, Theorem~\ref{thm:03.04} guarantees that the sets $\mcT_{j+2}$ and $\mcP_{j+2}$ are proper $C$--sets in $\R^{n+m}$ and $\R^n$. Since  $\mcC$ is a proper $C$--set in $\R^{n+m}$ and $\mcP_0^*$ is a $C$--set in $\R^n$, $\mcT_1^*=\mcC^*\oplus (A\ B)^T\mcP_0^*$ is a proper $C$--set in $\R^{n+m}$. In turn,  $\mcT_1$  and $\mcP_1=P_x\mcT_1$ are proper $C$--sets in $\R^{n+m}$ and $\R^n$. Hence, for all $k\in\N$, the sets $\mcT_{k+1}$ and $\mcP_{k+1}$ are proper $C$--sets in $\R^{n+m}$ and $\R^n$.

\subsection*{APPENDIX B--4: Proof of Proposition~\ref{prop:04.03}}
$(i)$: If for some $j\in\N$, $\mcP_{j+2}\subseteq \mcP_{j+1}$ or, equivalently,  $\mcP_{j+1}^*\subseteq \mcP_{j+2}^*$, then $\mcC^*\oplus (A\ B)^T\mcP_{j+1}^*\subseteq \mcC^*\oplus (A\ B)^T\mcP_{j+2}^*$ since $\mcC^*$ is a proper $C$--set in $\R^{n+m}$. It follows that $(\mcC^*\oplus (A\ B)^T\mcP_{j+2}^*)^*\subseteq (\mcC^*\oplus (A\ B)^T\mcP_{j+1}^*)^*$. In turn, $\mcP_{j+3}\subseteq \mcP_{j+2}$.  By Assumption $\mcP_0^*\subseteq \mcP_1^*$ so that $\mcP_2\subseteq \mcP_1$. Hence, for all $k\in \N$, $\mcP_{k+2}\subseteq \mcP_{k+1}$.\\
$(ii)$: The proof of this claim is conceptually identical to the proof of $(i)$. Note that assumption $\mcP_1^*\subseteq \mcP_0^*$ leads to reversed inclusions so that, for all $k\in \N$, $\mcP_{k+1}\subseteq \mcP_{k+2}$.

\subsection*{APPENDIX B--5: Proof of Proposition~\ref{prop:04.04}}
If for some $j\in\N$, $\mcP_{j+1}'\subseteq \mcP_{j+1}''$ or, equivalently,  $(\mcP_{j+1}'')^*\subseteq (\mcP_{j+1}')^*$, then $\mcC^*\oplus (A\ B)^T(\mcP_{j+1}'')^*\subseteq \mcC^*\oplus (A\ B)^T(\mcP_{j+1}')^*$ since $\mcC^*$ is a proper $C$--set in $\R^{n+m}$. In turn, $(\mcC^*\oplus (A\ B)^T(\mcP_{j+1}')^*)^*\subseteq (\mcC^*\oplus (A\ B)^T(\mcP_{j+1}'')^*)^*$. Thus, $\mcP_{j+2}'\subseteq \mcP_{j+2}''$.  By Assumption $(\mcP_0'')^*\subseteq (\mcP_0')^*$ so that $\mcP_1'\subseteq \mcP_1''$. Hence, for all $k\in \N$, $\mcP_{k+1}'\subseteq \mcP_{k+1}''$.

\subsection*{APPENDIX B--6: Proof of Proposition~\ref{prop:04.05}}
If for some $j\in \N$, $\mcL\subseteq \mcP_{j+1}\subseteq \mcD^+$, then, $(\mcD^+)^*\subseteq \mcP_{j+1}^*\subseteq \mcL^*$ and, since $\mcC^*$ is proper $C$--set in $\R^{n+m}$, $\mcC^*\oplus (A\ B)^T(\mcD^+)^*\subseteq \mcC^*\oplus (A\ B)^T\mcP_{j+1}^*\subseteq \mcC^*\oplus (A\ B)^T\mcL^*$. Since $\mcD^+=P_x\mcC$ is a proper $C$--set in $\R^n$, it follows that $\mcC^*\subseteq\mcC^*\oplus (A\ B)^T(\mcD^+)^*$. In turn, $\mcC^*\subseteq \mcC^*\oplus (A\ B)^T\mcP_{j+1}^*\subseteq \mcC^*\oplus (A\ B)^T\mcL^*$ and, consequently, $(\mcC^*\oplus (A\ B)^T\mcL^*)^*\subseteq (\mcC^*\oplus (A\ B)^T\mcP_{j+1}^*)^*\subseteq \mcC$. Hence,  $P_x(\mcC^*\oplus (A\ B)^T\mcL^*)^*\subseteq P_x(\mcC^*\oplus (A\ B)^T\mcP_{j+1}^*)^*\subseteq P_x\mcC$. By Proposition~\ref{prop:04.01}, $\mcL\subseteq P_x(\mcC^*\oplus (A\ B)^T\mcL^*)^*$. Also, $\mcD^+=P_x\mcC$. Hence, $\mcL\subseteq \mcP_{j+2}\subseteq \mcD^+$. By Assumption $\mcP_0^*\subseteq \mcL^*$, which, in turn, yields $\mcC^*\subseteq \mcC^*\oplus (A\ B)^T\mcP_0^*\subseteq \mcC^*\oplus (A\ B)^T\mcL^*$. Thus, $(\mcC^*\oplus (A\ B)^T\mcL^*)^*\subseteq (\mcC^*\oplus (A\ B)^T\mcP_0^*)^*\subseteq \mcC$. Hence, similarly as above, $P_x(\mcC^*\oplus (A\ B)^T\mcL^*)^*\subseteq P_x(\mcC^*\oplus (A\ B)^T\mcP_0^*)^*\subseteq P_x\mcC$. Since, $\mcL\subseteq P_x(\mcC^*\oplus (A\ B)^T\mcL^*)^*$ and $\mcD^+=P_x\mcC$, it follows that that $\mcL\subseteq \mcP_1\subseteq \mcD^+$. Hence, for all $k\in \N$, $\mcL\subseteq \mcP_{k+1}\subseteq \mcD^+$.

\subsection*{APPENDIX B--7: Proof of Theorem~\ref{thm:04.01}}
The space of nonempty compact subsets of $\R^n$ endowed with the Hausdorff distance is a complete metric space (see~\cite[Theorem~1.8.2.]{schneider:1993}). The space of nonempty convex compact subsets of $\R^n$ is a closed subset of the space of nonempty compact subsets of $\R^n$ (see~\cite[Theorem~1.8.5.]{schneider:1993}). The Blaschke selection theorem (see~\cite[Theorem~1.8.6.]{schneider:1993}) guarantees that any convergent sequence of nonempty convex compact subsets of $\R^n$ converges to a nonempty convex compact subset of $\R^n$.\\
$(i)$: By Propositions~\ref{prop:04.02},~\ref{prop:04.03}, and~\ref{prop:04.05} the sequence of the generator sets $\mcP_k$ is a sequence of inner and outer bounded and monotonic proper $C$--sets in $\R^n$. In particular, if $\mcP_0^*\subseteq \mcP_1^*$, then, for all $k\in \N$, $\mcL\subseteq \mcP_{k+2}\subseteq \mcP_{k+1}\subseteq \mcD^+$. Likewise,  if $\mcP_1^*\subseteq \mcP_0^*$, then, for all $k\in \N$, $\mcL\subseteq \mcP_{k+1}\subseteq \mcP_{k+2}\subseteq \mcD^+$. In either case, the Blaschke selection theorem, guarantees that  the sequence of the generator sets $\mcP_k$ converges to a nonempty compact convex subset $\mcP_\infty$ of $\R^n$. Since, for all $k\in \N$, $\mcL\subseteq \mcP_{k+1}\subseteq \mcD^+$, in the limit $\mcL\subseteq \mcP_\infty \subseteq \mcD^+$ so that $\mcP_\infty$ contains a proper $C$--set $\mcL$. In turn, the limit $\mcP_\infty$ is, in fact, a proper $C$--set in $\R^n$.\\
$(ii)$ For all $k\in \N$, $\mcT_{k+1}=\left(\mcC^*\oplus (A\ B)^T\mcP_k^*\right)^*$ are proper $C$--sets in $\R^{n+m}$ by Proposition~\ref{prop:04.02}. By~$(i)$ above, the sequence $\{\mcP_{k+1}\}_{k\in \N}$ of the generator sets $\mcP_k$ converges to the limit $\mcP_\infty$. In turn, the sequence $\{\mcP_{k+1}^*\}_{k\in\N}$ of polar sets $\mcP_k^*$ (each of which is a proper $C$--set in $\R^n$)  of generator sets $\mcP_k$ converges to the limit $\mcP_\infty^*$, which is a proper $C$--set in $\R^n$. The limit $\mcP_\infty^*$ is a proper $C$--set since for all $k\in \N$, $\mcL\subseteq \mcP_{k+1}\subseteq \mcD^+$ or, equivalently, $(\mcD^+)^*\subseteq \mcP_k^*\subseteq \mcL^*$ so that in the limit  $(\mcD^+)^*\subseteq \mcP_\infty^*\subseteq \mcL^*$. Thus, the limit $\mcP_\infty^*$ contains proper $C$--set $(\mcD^+)^*$. Consequently, the sequence $\{\mcC^*\oplus (A\ B)^T\mcP_k^*\}_{k\in \N}$ of proper $C$--sets in $\R^{n+m}$ converges to the limit $\mcC^*\oplus (A\ B)^T\mcP_\infty^*$ which is a proper $C$--set in $\R^{n+m}$. Clearly, the sequence $\{\mcT_{k+1}^*\}_{k\in \N}$ of polar sets $\mcT_k^*$ of the sets $\mcT_k$ is, by construction, such that, for all $k\in \N$, $\mcT_{k+1}^*=\mcC^*\oplus (A\ B)^T\mcP_k^*$. Thus, the sequence $\{\mcT_{k+1}^*\}_{k\in \N}$ of polar sets $\mcT_k^*$ of sets $\mcT_k$ converges to the proper $C$--set $\mcT_\infty^*=\mcC^*\oplus (A\ B)^T\mcP_\infty^*$ in $\R^{n+m}$. Thus, the sequence $\{\mcT_{k+1}\}_{k\ge 0}$  of proper $C$--sets $\mcT_k$ converges to a proper $C$--set $\mcT_\infty=(\mcT_\infty^*)^*$ in $\R^{n+m}$, and
\begin{equation*}
\mcT_\infty=\left(\mcC^*\oplus (A\ B)^T\mcP_\infty^*\right)^*.
\end{equation*}
Since the sequence $\{\mcT_{k+1}\}_{k\in\N}$ of sets $\mcT_k$  converges to the limit $\mcT_\infty$, which is a proper $C$--set in $\R^{n+m}$, the sequence $\{P_x\mcT_{k+1}\}_{k\in\N}$ of proper $C$--sets in $\R^n$ converges to the limit $P_x\mcT_\infty$, which is a proper $C$--set in $\R^n$. But, for all $k\in \N$, $\mcP_{k+1}=P_x\mcT_{k+1}$ so that  the sequence $\{\mcP_{k+1}\}_{k\in\N}$ of the generator sets $\mcP_k$ converges to this same limit $P_x\mcT_\infty$.  Hence, by the uniqueness of the limits,
\begin{equation*}
\mcP_\infty=P_x\mcT_\infty.
\end{equation*}
$(iii)$: By $(ii)$, $\mcP_\infty=P_x\mcT_\infty$ and $\mcT_\infty=\left(\mcC^*\oplus (A\ B)^T\mcP_\infty^*\right)^*$ so that
\begin{equation*}
\mcP_\infty=P_x\left(\mcC^*\oplus (A\ B)^T\mcP_\infty^*\right)^*,
\end{equation*}
and $\mcP_\infty$ solves the fixed point set--equation~\eqref{eq:04.03}.

\subsection*{APPENDIX B--8: Proof of Theorem~\ref{thm:04.02}}
Proposition~\ref{prop:04.04} guarantees that $\mcP_{k+1}'\subseteq \mcP_{k+1}''$ for all $k\in \N$. Theorem~\ref{thm:04.01} guarantees that the sequences of generator sets $\mcP_{k+1}'$ and $\mcP_{k+1}''$ converge to proper $C$--sets $\mcP_\infty'$ and $\mcP_\infty''$. Thus, as claimed, $\mcP_\infty'\subseteq \mcP_\infty''$.

\subsection*{APPENDIX C: Convergence of Value Functions $V_k\bb$}

When $\mcP_0^*$ is a $C$--set in $\R^n$ such that $\mcP_0^*\subseteq \mcP_1^*$, the sequence of generator sets $\mcP_k$ converges to $\mcP_\infty$. The sets $\mcP_k$ and the limit $\mcP_\infty$ are proper $C$--sets such that $\mcP_\infty\subseteq \mcP_{k+2}\subseteq \mcP_{k+1}$ for all $k\in\N$. Hence, for all $\delta >0$ there exists a  $k^*(\delta)\in \N$ such that, for all $k\ge k^*(\delta)$,
\begin{equation*}
\mcP_\infty\subseteq\mcP_k\subseteq (1+\delta)\mcP_\infty.
\end{equation*}
In turn, for all $k\ge k^*(\delta)$ and all $x\in \R^n$, $\gauge{\mcP_k}{x}\le \gauge{\mcP_\infty}{x}$ and $\gauge{(1+\delta)\mcP_\infty}{x}\le \gauge{\mcP_k}{x}$. The latter inequality states equivalently  $\gauge{\mcP_\infty}{x}\le (1+\delta)\gauge{\mcP_k}{x}$ since $\gauge{(1+\delta)\mcP_\infty}{x}=(1+\delta)^{-1}\gauge{\mcP_\infty}{x}$. Hence,  for all $k\ge k^*(\delta)$ and all $x\in \R^n$,
\begin{equation*}
0\le\gauge{\mcP_\infty}{x}-\gauge{\mcP_k}{x}\le \delta \gauge{\mcP_k}{x}\le \delta \gauge{\mcP_\infty}{x}.
\end{equation*}
Let $\phi_\infty:=\max_{x\in \mathbb{S}^{n-1}} \gauge{\mcP_\infty}{x}$. Note that $0<\phi_\infty <\infty$. 
Hence, for all $k\ge k^*(\delta)$ and all $x\in \mathbb{S}^{n-1}$,
\begin{equation*}
0\le\gauge{\mcP_\infty}{x}-\gauge{\mcP_k}{x}\le \delta \phi_\infty.
\end{equation*}
For all $\varepsilon>0$, $\delta(\varepsilon):=\varepsilon \phi_\infty^{-1}>0$ yields $\delta(\varepsilon) \phi_\infty=\varepsilon$. Hence, for all $\varepsilon>0$, there exists $k(\varepsilon)= k^*(\delta(\varepsilon))$ such that, for all $k\ge k(\varepsilon)$ and all $x\in \mathbb{S}^{n-1}$,
\begin{equation*}
0\le\gauge{\mcP_\infty}{x}-\gauge{\mcP_k}{x}\le \varepsilon.
\end{equation*}
Thus, the sequence of value functions $V_k\bb$ converges uniformly over the unit sphere $\mathbb{S}^{n-1}$ to $\gaugebb{\mcP_\infty}$.

When $\mcP_0^*$ is a $C$--set in $\R^n$ such that $\mcP_1^*\subseteq \mcP_0^*$, the argument is conceptually identical and, hence, omitted.

\subsection*{APPENDIX D: Convergence of Optimizer Maps $u_k\bb$ and Selections $\nu_k\bb$}

When $\mcP_0^*$ is a $C$--set in $\R^n$ such that $\mcP_0^*\subseteq \mcP_1^*$, the sequences of  sets $\mcP_k$ and $\mcT_k$ converge to $\mcP_\infty$ and $\mcT_\infty$. The sets $\mcP_k$ and $\mcT_k$ and the limits $\mcP_\infty$ and $\mcT_\infty$ are proper $C$--sets such that $\mcP_\infty\subseteq \mcP_{k+2}\subseteq \mcP_{k+1}$ and $\mcT_\infty\subseteq \mcT_{k+2}\subseteq \mcT_{k+1}$ for all $k\in\N$. Hence, for all $\delta >0$ there exists a  $k^*(\delta)\in \N$ such that for all $k\ge k^*(\delta)$
\begin{equation*}
\mcP_\infty\subseteq\mcP_k\subseteq (1+\delta)\mcP_\infty\text{ and }\mcT_\infty\subseteq\mcT_k\subseteq (1+\delta)\mcT_\infty.
\end{equation*}
Let $\phi_\infty:=\max_{x\in \mathbb{S}^{n-1}} \gauge{\mcP_\infty}{x}$. Note that $0<\phi_\infty <\infty$.  As shown in the Appendix~C, for all $k\ge k^*(\delta)$ and all $x\in \mathbb{S}^{n-1}$,
\begin{equation*}
\gauge{\mcP_k}{x}\le\gauge{\mcP_\infty}{x}\le \gauge{\mcP_k}{x}+ \delta \phi_\infty.
\end{equation*}
Thus, for all $k\ge k^*(\delta)$ and all $x\in \mathbb{S}^{n-1}$,
\begin{align*}
&\gauge{\mcP_k}{x}\mcT_k\subseteq \gauge{\mcP_\infty}{x}\mcT_\infty\oplus \delta\gauge{\mcP_\infty}{x}\mcT_\infty\text{ and}\\
&\gauge{\mcP_\infty}{x}\mcT_\infty\subseteq \gauge{\mcP_k}{x}\mcT_k\oplus \delta\phi_\infty\mcT_k.
\end{align*}
Let, for all $k\in \N$,  $\tau_{k+1}:=H_{\mcB^{n+m}}(\mcT_{k+1},\{0\})$ and $\tau_\infty:=H_{\mcB^{n+m}}(\mcT_\infty,\{0\})$. Note that $0<\tau_\infty\le\tau_{k+2}\le\tau_{k+1}<\infty$  for all $k\in\N$ since,  for all $k\in\N$, $\mcT_\infty\subseteq \mcT_{k+2}\subseteq \mcT_{k+1}$ and $\mcT_\infty$ and $\mcT_{k+1}$ are proper $C$--sets. It follows that,  for all $k\ge k^*(\delta)$ and all $x\in \mathbb{S}^{n-1}$,
\begin{align*}
&\gauge{\mcP_k}{x}\mcT_k\subseteq \gauge{\mcP_\infty}{x}\mcT_\infty\oplus \delta \phi_\infty\tau_1\mcB^{n+m}\text{ and}\\
&\gauge{\mcP_\infty}{x}\mcT_\infty\subseteq \gauge{\mcP_k}{x}\mcT_k\oplus \delta \phi_\infty\tau_1\mcB^{n+m}.
\end{align*}
For all $\varepsilon>0$, $\delta(\varepsilon):=\varepsilon \phi_\infty^{-1}\tau_1^{-1}>0$ yields $\delta(\varepsilon) \phi_\infty \tau_1=\varepsilon$. Hence, for all $\varepsilon>0$, there exists $k(\varepsilon)= k^*(\delta(\varepsilon))$ such that, for all $k\ge k(\varepsilon)$ and all $x\in \mathbb{S}^{n-1}$,
\begin{align*}
&\gauge{\mcP_k}{x}\mcT_k\subseteq \gauge{\mcP_\infty}{x}\mcT_\infty\oplus \varepsilon\mcB^{n+m}\text{ and}\\
&\gauge{\mcP_\infty}{x}\mcT_\infty\subseteq \gauge{\mcP_k}{x}\mcT_k\oplus \varepsilon\mcB^{n+m},
\end{align*}
and the sequence of set--valued maps $\gaugebb{\mcP_k}\mcT_k$ converges uniformly over the unit sphere $\mathbb{S}^{n+m-1}$ to set--valued map $\gaugebb{\mcP_\infty}\mcT_\infty$. 

Hence, since $u_k(x)=\{u\in \R^m\ :\ (x,u)\in\gauge{\mcP_k}{x}\mcT_k\}$ as specified in~\eqref{eq:04.10} and $u_\infty(x)=\{u\in \R^m\ :\ (x,u)\in\gauge{\mcP_\infty}{x}\mcT_\infty\}$ as specified in~\eqref{eq:04.21}, the sequence of optimizer maps $u_k\bb$ converges uniformly over the unit sphere $\mathbb{S}^{n-1}$ to $u_\infty\bb$.

Likewise, since $\nu_k(x)=\arg\min_u\{u^Tu\ :\ (x,u)\in\gauge{\mcP_k}{x}\mcT_k\}$ as specified in~\eqref{eq:04.10} and $\nu_\infty(x)=\arg\min_u\{u^Tu\ :\ (x,u)\in\gauge{\mcP_\infty}{x}\mcT_\infty\}$ as specified in~\eqref{eq:04.21}, the sequence of selections $\nu_k\bb$ converges uniformly over the unit sphere $\mathbb{S}^{n-1}$ to $\nu_\infty\bb$.

When $\mcP_0^*$ is a $C$--set in $\R^n$ such that $\mcP_1^*\subseteq \mcP_0^*$, the argument is conceptually identical and, hence, omitted.

\bibliographystyle{unsrt}
\bibliography{MBE}
\end{document}